\documentclass[11pt]{amsart}
\theoremstyle{plain}

\begin{document}
\large

\title[Monotonicity Formulae]{New Monotonicity Formulae for Semi-linear
Elliptic and Parabolic Systems}

\author{Li Ma}
\address{Department of Mathematic Science, Tsinghua University,
 Peking 100084, P. R. China }

\email{lma@math.tsinghua.edu.cn}
\author{Xianfa Song }
\address{Department of Mathematic Science, Tsinghua University,
 Peking 100084, P. R. China }
\email{xsong@math.tsinghua.edu.cn}

\author{Lin Zhao}
\address{Department of Mathematic Science, Tsinghua University,
 Peking 100084, P. R. China }

\begin{abstract}
In this paper, we establish a general monotonicity formula of the
following elliptic system
$$
\Delta u_i+f_i(u_1,...,u_m)=0 \quad {\rm in} \ \Omega, \label{0.1}
$$
where $\Omega\subset\subset \mathbb{R}^n$ is a bounded domain,
$(f_i(u_1,...,u_m))=\nabla F(\vec{u})$, and $F(\vec{u})$ is a given
smooth function of $\vec{u}=(u_1,...,u_m)$, $m,n$ are two positive
integers. We also set up a new monotonicity formula for the
following parabolic system
$$
\partial_t u_i-\Delta u_i-f_i(u_1,...,u_m)=0 \quad {\rm in} \ (t_1, t_2)\times
\mathbb{R}^n,
$$
where $t_1<t_2$ are two constants, $(f_i(\vec{u}))$ is given as
above. Our new monotonicity formulae are focused on more attention
to the monotonicity of non-linear terms. Our point of view is that
we introduce an index called $\beta$ to measure the monotonicity of
the non-linear terms in the problems. The index in the study of
monotonicity formulae is very useful in understanding the behavior
of blow up sequences of solutions. Corresponding monotonicity
results for free boundary problems are also presented.

\end{abstract}

\thanks{This work is partially
supported by the key project 973 of the Ministry of Sciences and
Technology of China and partially supported by the grant from
China Postdoctoral Science Foundation. The fax number of Li Ma is
8610-62781785}

\maketitle
\date{}

\newtheorem{theorem}{Theorem}[section]
\newtheorem{definition}{Definition}[section]
\newtheorem{lemma}{Lemma}[section]
\newtheorem{proposition}{Proposition}[section]
\newtheorem{corollary}{Corollary}[section]
\newtheorem{remark}{Remark}[section]
\renewcommand{\theequation}{\thesection.\arabic{equation}}
\catcode`@=11 \@addtoreset{equation}{section} \catcode`@=12

\section{introduction}
In this paper, we will establish a general monotonicity formula of
the following elliptic system
\begin{align}
\Delta u_i+f_i(u_1,...,u_m)=0 \quad {\rm in} \ \Omega, \label{0.1}
\end{align}
where $\Omega\subset\subset \mathbb{R}^n$ is a bounded domain,
$(f_i(u_1,...,u_m))=\nabla F(\vec{u})$, and $F(\vec{u})$ is a given
smooth function of $\vec{u}=(u_1,...,u_m)$, $m,n$ are two positive
integers. Here we assume that the solution $\vec{u}\in
H^1_{loc}(\Omega)$ satisfies (\ref{0.1}) in the variational sense to
be defined in section two. We remark that smooth solutions to
(\ref{0.1}) satisfy (\ref{0.1}) in the variational sense. We will
also establish a monotonicity formula for regular solutions of the
following parabolic system
\begin{align}
\partial_t u_i-\Delta u_i-f_i(u_1,...,u_m)=0 \quad {\rm in} \ (t_1, t_2)\times
\mathbb{R}^n, \label{0.2}
\end{align}
where $t_1<t_2$ are two constants, $(f_i(\vec{u}))$ is given as
above. We also consider corresponding results for free boundary
problems. The new point in our monotonicity formula is that we
introduce an index $\beta$, which measures the monotonicity of the
non-linear term $f_i$. This index $\beta$ also gives us the rate of
scaled sequence of the blow-up process for implied solutions. Our
main results are Theorem 2.1, Theorem 2.2, Theorem 3.1, and Theorem
3.3 below. As a corollary, we can also give a monotonicity formula
for Ginzburg-Landau model (see our Assertion 2.1 below). Here we
give a brief introduction to our results.

Before we state the monotonicity formulae, we introduce some
notations and concepts. We will use some notations of \cite{W3} in
convenience, and denote by $x\cdot y$ the Euclidean inner product in
$\mathbb{R}^n\times \mathbb{R}^n$, by $|x|$ the Euclidean norm in
$\mathbb{R}^n$, by $B_r(x_0):=\{x\in \mathbb{R}^n| |x-x_0|<r\}$ the
ball of center $x_0$ and radius $r$, by
$Q_r(x_0,t_0):=(t_0-r^2,t_0+r^2)\times B_r(x_0)$ the cylinder of
radius $r$ and height $2r^2$, by
$T^-_r(t_0):=(t_0-4r^2,t_0-r^2)\times \mathbb{R}^n$ the horizontal
layer from $t_0-4r^2$ to $t_0-r^2$, by $T^+_r(t_0):=(t_0+r^2,
t_0+4r^2)\times \mathbb{R}^n$ the horizontal layer from $t_0+r^2$ to
$t_0+4r^2$, and by
$$G_{(t_0,x_0)}(t,x):=4\pi(t_0-t)|4\pi(t_0-t)|^{-\frac{n}{2}-1}
\exp\left(-\frac{|x-x_0|^2}{4(t_0-t)}\right)
$$
the backward heat kernel, defined in $((-\infty,
t_0)\cup(t_0,+\infty))\times \mathbb{R}^n$. Sometimes, we denote by
$T^-_r:=(T-4r^2,T-r^2)\times \mathbb{R}^n$ and $T^+_r:=(T+r^2,
T+4r^2)\times \mathbb{R}^n$ for $t_0=T$. Furthermore, by $\nu$ we
will always refer to the outer unit normal on a given surface. We
denote by $\mathcal{H}^s$ the $s$-dimensional Hausdorff measure, and
$H^1_{loc}(\Omega)$ and $H^1(Q_T)$ the usual Sobolev space and
parabolic Sobolev spaces respectively as defined in \cite{LSU}.

Roughly speaking, our new monotonicity formula for (\ref{0.1}) is as
follows. We will show that for the variational solution $\vec{u}$ to
(\ref{0.1}), the function
$$
\Phi_{x_0}(r):=r^{-n-2\beta+2}\int_{B_r(x_0)}(|\nabla
\vec{u}|^2-2F(\vec{u}))-\beta r^{-n-2\beta+1}\int_{\partial
B_r(x_0)} \vec{u}^2d\mathcal{H}^{n-1}
$$
is increasing in $r$ if
\begin{equation}
\int_{B_r(x_0)}[2(\beta-1)F(\vec{u})-\beta
\vec{u}\vec{f}(\vec{u})]\geq 0;\label{C1}
\end{equation}
and for (\ref{0.2}) the functions
\begin{align*}
\Psi^-(r)&=r^{-2\beta}\int_{T^-_r}(|\nabla
\vec{u}|^2-2F(\vec{u}))G_{(T,x_0)}
-\frac{\beta}{2}r^{-2\beta}\int_{T^-_r}\frac{1}{T-t}\vec{u}^2G_{(T,x_0)}
\end{align*}
and
\begin{align*}
\Psi^+(r)&=r^{-2\beta}\int_{T^+_r}(|\nabla
\vec{u}|^2-2F(\vec{u}))G_{(T,x_0)}
-\frac{\beta}{2}r^{-2\beta}\int_{T^+_r}\frac{1}{T-t}\vec{u}^2G_{(T,x_0)}
\end{align*}
are increasing in $r$ for $\beta$ such that
\begin{equation}
\int_{T^-_r}[2(\beta-1)F(\vec{u})-\beta
\vec{u}\vec{f}(\vec{u})]G_{(T,x_0)}\geq 0 \label{C2}
\end{equation}
and
\begin{equation}
\int_{T^+_r}[2(\beta-1)F(\vec{u})-\beta
\vec{u}\vec{f}(\vec{u})]G_{(T,x_0)}\geq 0, \label{C3}
\end{equation}
where $\vec{f}=(f_1,...,f_m)$.

We remark that conditions (\ref{C1})(\ref{C2})(\ref{C3}) are
automatically true if
$$ 2(\beta-1)F(\vec{u})-\beta
\vec{u}\vec{f}(\vec{u})\geq 0,\quad {\rm for} \ \vec{u}\in
\mathbb{R}^n.
$$
We will give more illustration by examples in section 2 and section
4. From the expression above, it is clear that the number $\beta$
measures the monotonicity of the non-linear term, and our new
monotonicity formulae are focused on more attention to the
monotonicity of non-linear terms. We emphasize that the boundary
term in the elliptic case is important, and in some special cases,
it was noticed by Weiss (see \cite {W}, \cite{W1}, \cite{W2} and
\cite{W3}) who called it ``boundary-adjusted energy''. This term is
nature in measuring the flux transportation through boundary. Our
method can also be used to study elliptic /parabolic systems with
variable coefficients. For example, one may extend the monotonicity
results above to elliptic systems and parabolic systems with
variable coefficients.

As corollaries of our results to systems (\ref{0.1}) and
(\ref{0.2}), we subsequently establish the monotonicity formulae for
the following general elliptic equation
\begin{align}
 \Delta u+f(u)=0 \quad {\rm in} \ \Omega, \label{0.3}
\end{align}
where $\Omega\subset\subset \mathbb{R}^n$ is a bounded domain, and
the parabolic equation
\begin{align}
\partial_t u-\Delta u=f(u), \quad {\rm in} \ (t_1, t_2)\times
\mathbb{R}^n, \label{0.4}
\end{align}
where $t_1,t_2$ are two constants, $f(u)$ is a given function of
$u$. The results will be stated in detail in section four.

As we said, our present work is closely related to the monotonicity
formula of Weiss \cite {W}, \cite{W1}, \cite{W2} and \cite{W3} and
the monotonicity formula of Alt, Caffarelli and Friedman \cite{ACF}.
However, we will not only obtain the monotonicity formulae for more
general models of single equation, but also establish the
monotonicity formulae for some types of elliptic and parabolic
systems, and the results are completely new. Moreover, we can choose
different $\beta$ such that the monotonicity formula holds even in
the same model of Weiss' papers (\cite{W1},\cite{W2},\cite{W3}). For
example, if $f(u)=u^p$ with $-1<p<1$, we can choose any $\beta\geq
2/(1-p)$; while for $p<-1$ or $p>1$, we can choose any $\beta\leq
2/(1-p)$. And we can construct different types of the scaled
sequences through the choosing of $\beta$. For example, denoting the
sequences $u_k(x):=\rho_k^{-\beta}u(x+\rho_k x)$ for $\beta<0$, we
find that they are different from the {\it blow up sequences} for
$\beta>0$.

We now further compare our result with those of Weiss and give a
brief review about monotonicity formulae related. In \cite {W},
\cite{W1}, \cite{W2} and \cite{W3}, Weiss introduced the
``boundary-adjusted energy'', and obtained some new monotonicity
formulae. In \cite{W1}, Weiss studied the critical points with
respect to the energy
$$
w \rightarrow F(w)=\int_{\Omega}(|\nabla
w|^2+\lambda_+\chi_{\{w>0\}}w^p+\lambda_-\chi_{\{w<0\}}(-w)^p)
$$
with $p\in [0,2)$ and found that: Assume that $u$ is a solution and
$B_{\delta}(x_0)\subset \Omega$. Then, in term of our results, for
$\beta:=\frac{2}{2-p}$ and for any $0<\rho<\sigma<\delta$ the
function
\begin{align*}
\Phi(r)&:=r^{-n-2\beta+2}\int_{B_r(x_0)}(|\nabla
u|^2+\lambda_+\chi_{\{u>0\}}u^p+\lambda_-\chi_{\{u<0\}}(-u)^p)\\
&\qquad -\beta r^{-n-2\beta+1}\int_{\partial B_r(x_0)}
u^2d\mathcal{H}^{n-1}
\end{align*}
defined in $(0,\delta)$, satisfies the monotonicity formula
$$
\Phi(\sigma)-\Phi(\rho)=\int^{\sigma}_{\rho}
r^{-n-2\beta+2}\int_{\partial B_r(x_0)}2(\nabla u\cdot
\nu-\beta\frac{u}{r})^2d\mathcal{H}^{n-1}dr\geq 0.$$ In \cite{W2},
the monotonicity formula for $\Delta u=\chi_{\{u>0\}}$ has the same
form of $\Phi(r)$ with $p=1$.

In \cite{W3}, Weiss studied the gradient flow in
$L^2(\mathbb{R}^n)$ with respect to the energy
$$
w \rightarrow F(w)=\int_{\mathbb{R}^n}(|\nabla
w|^2+\lambda_+\chi_{\{w>0\}}w^p+\lambda_-\chi_{\{w<0\}}(-w)^p)
$$
with $p\in [0,2)$ and found that: Assume that $t_1\leq T\leq t_2$,
$x_0\in \mathbb{R}^n$ and $u$ is a solution with some conditions.
Then, again in terms of our results, for $\beta:=\frac{2}{2-p}$ and
for any $0<\rho<\sigma<\delta$ the function
\begin{align*}
\Psi^-(r)&:=r^{-2\beta}\int_{T^-_r(T)}(|\nabla
u|^2+\lambda_+\chi_{\{u>0\}}u^p+\lambda_-\chi_{\{u<0\}}(-u)^p)G_{(T,x_0)}\\
&\qquad -\frac{\beta}{2} r^{-2\beta}\int_{T^-_r(T)}
\frac{1}{T-t}u^2G_{(T,x_0)}
\end{align*}
and
\begin{align*}
\Psi^+(r)&:=r^{-2\beta}\int_{T^+_r(T)}(|\nabla
u|^2+\lambda_+\chi_{\{u>0\}}u^p+\lambda_-\chi_{\{u<0\}}(-u)^p)G_{(T,x_0)}\\
&\qquad -\frac{\beta}{2} r^{-2\beta}\int_{T^+_r(T)}
\frac{1}{T-t}u^2G_{(T,x_0)}
\end{align*}
are well defined in the interval $(0,\frac{\sqrt{T-t_1}}{2})$ and
$(0,\frac{\sqrt{t_2-T}}{2})$, respectively, and satisfy for any
$0<\rho<\sigma<\frac{\sqrt{T-t_1}}{2}$ and
$0<\rho<\sigma<\frac{\sqrt{t_2-T}}{2}$, respectively, the
monotonicity formulae
\begin{align*}
&\quad \Psi^-(\sigma)-\Psi^-(\rho)\\
&=\int^{\sigma}_{\rho}r^{-2\beta-1}
\int_{T^-_r(T)}\frac{1}{T-t}(\nabla
u\cdot(x-x_0)-2(T-t)\partial_t u-\beta u)^2G_{(T,x_0)}\\
&\geq 0
\end{align*}
and
\begin{align*}
&\quad
\Psi^+(\sigma)-\Psi^+(\rho)\\
&=\int^{\sigma}_{\rho}r^{-2\beta-1}
\int_{T^+_r(T)}\frac{1}{T-t}(\nabla
u\cdot(x-x_0)-2(T-t)\partial_t u-\beta u)^2G_{(T,x_0)}\\
&\geq 0.
\end{align*}

 In \cite{ACF}, Alt, Caffarelli and Friedman established a
monotonicity formula for variational problems with two phases and
their free boundaries. The monotonicity formula of Alt-Caffarelli
-Friedman
 plays an
important role as a fundamental and powerful tool in free boundary
problems. Roughly speaking, they found that
$$
\Phi(r)=\left(\frac{1}{r^2}\int_{B_r(x_0)}\frac{|\nabla
h_1|^2}{|x-x_0|^{N-2}}\right)\left(\frac{1}{r^2}\int_{B_r(x_0)}\frac{|\nabla
h_2|^2}{|x-x_0|^{N-2}}\right)$$ is increasing in $r(0<r<R)$ for the
sub-solutions $h_1,h_2$ of $\Delta u=0$ in $B(x_0, R)(R>0)$ with
$h_1h_2=0$ and $h_1(x_0)=h_2(x_0)=0$. We can also see \cite{CKS}. In
\cite{CJK}, Caffarelli, Jerison and Kenig found that there is a
dimensional constant $C$ such that
\begin{align*}
\Phi(r)&=\left(\frac{1}{r^2}\int_{B_r}\frac{|\nabla
u_+|^2}{|X|^{n-2}}dX\right)\left(\frac{1}{r^2}\int_{B_r}\frac{|\nabla
u_-|^2}{|X|^{n-2}}dX\right)\\
&\leq C\left(1+\int_{B_1}\frac{|\nabla
u_+(X)|^2}{|X|^{n-2}}dX+\int_{B_1}\frac{|\nabla
u_-(X)|^2}{|X|^{n-2}}dX\right)^2
\end{align*}
with $ 0<r\leq 1$ for $\Delta u_{\pm}\geq -1$ in the sense of
distributions, and $u_+(X)u_-(X)=0$ for all $X\in B_1$.

Various monotonicity formulae have catched many authors' attentions
in the past several years. Let's us briefly review some progress in
them. The well-known monotonicity formula for minimal hyper-surfaces
in \cite{Si}
$$
\frac{d}{dr}\left(\frac{\mathcal{H}^n(M\cap B_r)}{r^n} \right)
=\frac{d}{dr}\int_{M\cap
B_r}\frac{|x^{\perp}|^2}{|x|^{n+2}}d\mathcal{H}^n$$ is a local
statement in balls $B_r\subset \mathbb{R}^{n+1}$, which plays a
very important role in analyzing singularity set. There are many
references about the topic.  Fleming obtained the monotonicity
formula for area minimizing currents in \cite{Fl}. Allard proved
the monotonicity formula for stationary rectifiable $n$-varifolds
in \cite{Al}. Schoen and Uhlenbeck established the monotonicity
formula for harmonic maps in \cite{SU}. Price proved the
monotonicity for weakly stationary harmonic maps and Yang-Mills
equations in \cite{Pr}. Giga and Kohn obtained in \cite{GK} the
monotonicity formula for the solutions of semi-linear heat
equations $\partial_t u-\Delta u-|u|^{p-1}u=0$ with blow-up
analysis, where $p>1$, and Pacard established its localization for
weakly stationary solutions of the corresponding elliptic equation
 in \cite{Pa}.  M.Struwe derived the monotonicity formula involving the
associated energy densities for the equation $\partial_t u-\Delta
u\in T^{\perp}N$ in \cite{St}. Riviere \cite{RI}, F.H.Lin and
Riviere \cite{LR}, Bourgain, Brezis, and Mironescu \cite{BBM} set up
some monotonicity formulae for Ginzburg-Landau model. The famous
monotonicity formula for mean curvature flow, which was found by
G.Huisken \cite{Hu}, says that
$$\frac{d}{dt}\int_{M_t}G
d\mu_t=-\int_{M_t}|\vec{H}-\frac{x^{\perp}}{2t}|^2G d\mu_t,$$
which involves the backward heat kernel function
$G(x,t)=\frac{1}{(-4\pi t)^{n/2}}e^{\frac{|x|^2}{4t}}$ for $t<0$
and $x\in \mathbb{R}^{n+k}$. Monotonicity formulae for geometric
evolution equations on more general domains were also derived by
Hamilton in \cite{Ha}. In \cite{E1} and \cite{E2}, the local
monotonicity formula had been given by Ecker in the "heat-ball"
$$
E^{\gamma}_r=\{(x,t)\in \mathbb{R}^n\times \mathbb{R}, t<0,
\Phi^{\gamma}>\frac{1}{r^{n-\gamma}}\}=\bigcup_{-\frac{r^2}{4\pi}<t<0}
B_{R^{\gamma}_r(t)}\times\{t\},
$$
where
$$ \Phi^{\gamma}(x,t)=\frac{1}{(-4\pi
t)^{\frac{n-\gamma}{2}}}e^{\frac{|x|^2}{4t}}, \quad
R^{\gamma}_r(t)=\sqrt{2(n-\gamma)\log(\frac{-4\pi t}{r^2})}.
$$
It can be written as follows:
\begin{align*} \frac{d}{dr} \{ \frac{1}{R^{n-\gamma
}} \int_{E^{\gamma }_r } & \frac{n-\gamma }{-2t}( e(u)-\frac{\beta
}{2t}u^2) -\frac{x}{2t} \cdot Du (\frac{\partial u}{\partial
t}+\frac{x}{2t} \cdot Du+\frac{\beta }{t} u)dxdt \} \\
&=\frac{n-\gamma }{r^{n-\gamma+1}} \int_{E^{\gamma }_r}
\left(\frac{\partial u}{\partial t} +\frac{x}{2t}\cdot
Du+\frac{\beta }{t}u \right)^2dxdt,
\end{align*}
where $u$ is a solution of $u_t-\Delta u-|u|^{p-1}u=0, \ x\in
\mathbb{R}^n, t<0 $ with $p>1$.

The monotonicity formula also appears in the parabolic potential
theory \cite{CPS}. For a function $v$ and any $t>0$, define
$$
I(t;v)=\int^0_{-t}\int_{\mathbb{R}^n} |\nabla
v(s,x)|^2G(-s,x)dxds.$$  In \cite{Ca}, Caffarelli found that
$$
\Phi(t)=\Phi(t;h_1,h_2):=\frac{1}{t^2}I(t;h_1)I(t;h_2)$$ is
monotone nondecreasing in $t(0<t<1)$ for  nonnegative subcaloric
functions $h_1, h_2$ in the strip $[-1,0]\times \mathbb{R}^n$, $
h_1(0,0)=h_2(0,0)=0$ and $ h_1\cdot h_2=0$ with a polynomial
growth at infinity. Its localization can be stated as follows:
There exists a constant $C=C(n,\psi)>0$ such that
$$
\Phi(t;w_1,w_2)\leq C\|h_1\|^2_{L^2(Q^-_1)}\|h_2\|^2_{L^2(Q^-_1)}
$$
for any $0<t<1/2$, here $\psi(x)\geq 0$ be a $C^{\infty}$ cut-off
function with {\rm supp}$\psi\subset B_{3/4}$ and
$\psi|_{B_{1/2}}=1$ and $w_i=h_i\psi$, see \cite{CPS}. In \cite{CK},
this formula was generalized  for parabolic equations with variable
coefficients, and was written as
\begin{align*}
&\frac{1}{t}\int^0_{-t}\int_{\mathbb{R}^n}
|\nabla(u_1\psi)|^2G(x,-s)dxds\cdot
\frac{1}{t}\int^0_{-t}\int_{\mathbb{R}^n} |\nabla(u_2\psi)|^2G(x,-s)dxds\\
&\quad \leq C\left(\|u_1\|^4_{L^2(Q_2)}+\|u_2\|^4_{L^2(Q_2)}
\right).
\end{align*}

The remaining part of the paper is organized as follows. In section
2 we establish the monotonicity formula for (\ref{0.1}) and
characterize the scaled sequences. In section 3 we establish the
monotonicity formula for (\ref{0.2}) and characterize the scaled
sequences. In section 4 we state the monotonicity formulae for
(\ref{0.3}) and (\ref{0.4}) and give some examples.

\section{the monotonicity formula of an elliptic system}
Consider the elliptic system
\begin{align} \Delta u_i+f_i(u_1,..., u_m)=0, \
i=1,...,m, \quad {\rm in} \ \Omega,\label{2.1}
\end{align}
where $\Omega\subset \subset\mathbb{R}^n$ and $(f_i)$ is the
gradient of a given smooth function $F(u_1,..., u_m)$.

In order to define the \emph{variational solution} of (\ref{2.1}),
we need to give some notations. We denote by  for
$\phi=(\phi_1,...,\phi_n)\in
H^{1}_{loc}(\mathbb{R}^n;\mathbb{R}^n)$
$$
D\phi:=\left(\begin{array}{ll} \partial_1 \phi_1 \ ... \
&\partial_n \phi_1\\
\qquad \ ...& \\
\partial_1 \phi_n \ ... \ &\partial_n \phi_n
\end{array} \right)
$$
and by $\nabla \vec{u}\cdot x=(\nabla u_1\cdot x,...,\nabla
u_m\cdot x)$, and by $(\nabla \vec{u}\cdot
\nu)^2=\sum_{i=1}^m(\nabla u_i\cdot \nu)^2$, by $\vec{u}\nabla
\vec{u}\cdot \nu=\sum_{i=1}^m u_i\nabla u_i\cdot \nu$ for any
vector $\nu$, and by $\nabla \vec{u} D\phi \nabla
\vec{u}=\sum_{i=1}^m \nabla u_i D\phi \nabla u_i$. We say
$\vec{u}\in C^0(\Omega)\cap C^2(\Omega)$ if every $u_i\in
C^0(\Omega)\cap C^2(\Omega)$. We say $\vec{u}\in
H^{1}_{loc}(\Omega)$ if every $u_i\in H^{1}_{loc}(\Omega)$,
$i=1,... ,m$.

{\bf Definition 1.} {\it We call $\vec{u}\in H^{1}_{loc}(\Omega)$ is
a solution of (\ref{2.1}) in the sense of variations, or simply a
variational solution, if
 $\vec{u}\in H^1_{loc}(\Omega)$ satisfies (\ref{2.1}) in the distributional sense
 with
 $$u_if_i(\vec{u}),\;\;F(\vec{u})\in L^1_{loc}(\Omega)$$
 for $i=1,...,m$,
and the first variation with respect to domain variations of the
functional
$$
G(\vec{v}):=\int_{\Omega}(|\nabla \vec{v}|^2-2F(\vec{v}))$$
vanishes at $\vec{v}=\vec{u}$, i.e.
\begin{align*}
0=\frac{d}{d\varepsilon}G(\vec{u}(x+\varepsilon
\phi(x)))|_{\varepsilon=0}=\int_{\Omega} [(|\nabla
\vec{u}|^2-2F(\vec{u})) {\rm div}\phi-2\nabla \vec{u} D\phi \nabla
\vec{u}]
\end{align*}
for any $\phi \in C^1_0(\Omega; \mathbb{R}^n)$. }

\begin{theorem}\label{emf05}
 Assume that $\vec{u}$ is a solution of (\ref{2.1}) in the sense of variations
 in the ball
 $B_{\delta}(x_0)\subset\subset \Omega$. Then for any
 $\beta\in
\mathbb{R}$ such that
$$\int_{B_r(x_0)}(2(\beta-1)F(\vec{u})-\beta
\vec{u} \vec{f}(\vec{u}))\geq 0\quad  {\rm for} \ 0<r\leq \delta
$$
the function $$ \Phi_{x_0}(r):=r^{-n-2\beta+2}\int_{B_r(x_0)}
(|\nabla \vec{u}|^2-2F(\vec{u}))-\beta r^{-n-2\beta+1}\int_{\partial
B_r(x_0)} \vec{u}^2d\mathcal{H}^{n-1},
$$
defined in $(0,\delta)$, satisfies the monotonicity formula
\begin{align}
\Phi_{x_0}(\sigma)-\Phi_{x_0}(\rho)&=
\int_{\rho}^{\sigma}2r^{-n-2\beta+1}\int_{B_r(x_0)}(2(\beta-1)F(\vec{u})-\beta
\vec{u} \vec{f}(\vec{u}))\label{2.2}\\
&\ +\int_{\rho}^{\sigma}2r^{-n-2\beta+2}\int_{\partial B_r(x_0)}
(\nabla \vec{u}\cdot \nu-\beta \frac{\vec{u}}{r})^2
d\mathcal{H}^{n-1}dr\geq 0\nonumber
\end{align}
for all $ 0<\rho<\sigma<\delta$, where
$$(\nabla \vec{u}\cdot \nu-\beta
\frac{\vec{u}}{r})^2=\sum_{i=1}^m(\nabla u_i\cdot \nu-\beta
\frac{u_i}{r})^2.
$$
\end{theorem}

\begin{proof}
We may assume that $x_0=0$ by a translation. We take after
approximation $\phi_{\varepsilon}(x):=\eta_{\varepsilon}(x)x$ as
test function in Definition 1 for small positive $\varepsilon$ and
$\eta_{\varepsilon}(x):=\max(0,
\min(1,\frac{r-|x|}{\varepsilon}))$, and obtain that
\begin{align}
0&=\int[n(|\nabla \vec{u}|^2-2F(\vec{u}))\eta_{\varepsilon}
-2|\nabla \vec{u}|^2\eta_{\varepsilon}]\nonumber\\
&\qquad+\int(|\nabla \vec{u}|^2-2F(\vec{u}))\nabla
\eta_{\varepsilon}\cdot x-2\nabla \vec{u}\cdot x
\nabla \vec{u} \cdot \nabla \eta_{\varepsilon})\nonumber\\
&\rightarrow \int_{B_r(0)} [n(|\nabla \vec{u}|^2-2F(\vec{u}))-2|\nabla \vec{u}|^2]\nonumber\\
&\qquad -\int_{\partial B_r(0)}[r(|\nabla
\vec{u}|^2-2F(\vec{u}))-2r(\nabla \vec{u}\cdot
\nu)^2]d\mathcal{H}^{n-1}\label{2.3}
\end{align}
for a.e. $r\in (0,\delta)$ as $\varepsilon\rightarrow 0$.

Using mollifier $u_{i,\rho}$ to (\ref{2.1}) for every $u_i$
$(i=1,...,m)$, where $\rho>0$, we have
$$
-\Delta u_{i,\rho}=(f_i(u_1,...,u_m))_{\rho}.
$$
Multiplying this equation by $u_{i}$ and integrating over
$B_r(0)$, then sending $\rho\to 0+$ , we can easily derive the
formula
\begin{equation} \int_{B_r(0)}|\nabla \vec{u}|^2=\int_{\partial B_r(0)}
\vec{u}\nabla \vec{u}\cdot \nu
d\mathcal{H}^{n-1}+\int_{B_r(0)}\vec{u}\vec{f}(\vec{u})\label{2.4}
\end{equation}
for a.e. $r\in(0,\delta)$. Next, multiplying (\ref{2.3}) by
$-r^{-n-2\beta+1}$ and using (\ref{2.4}), we obtain that
\begin{align*}
0&=-r^{-n-2\beta+1}\int_{B_r(0)} [n(|\nabla
\vec{u}|^2-2F(\vec{u}))-2|\nabla \vec{u}|^2]\\& \quad
+r^{-n-2\beta+2}\int_{\partial B_r(0)}[(|\nabla
\vec{u}|^2-2F(\vec{u}))-2(\nabla \vec{u}\cdot \nu)^2]d\mathcal{H}^{n-1}\\
&=(-n-2\beta+2)r^{-n-2\beta+1}\int_{B_r(0)} (|\nabla
\vec{u}|^2-2F(\vec{u}))\\&\quad-2(2\beta-2)r^{-n-2\beta+1}\int_{B_r(0)}F(\vec{u})\\
&\quad+2\beta r^{-n-2\beta+1}(\int_{\partial B_r(0)} \vec{u}\nabla
\vec{u}\cdot \nu d\mathcal{H}^{n-1}+\int_{B_r(0)} \vec{u}\vec{f}(\vec{u}))\\
&\quad +r^{-n-2\beta+2}\int_{\partial B_r(0)}(|\nabla
\vec{u}|^2-2F(\vec{u}))d\mathcal{H}^{n-1}\\
&\quad -2r^{-n-2\beta+2}\int_{\partial B_r(0)} (\nabla \vec{u} \cdot
\nu)^2d\mathcal{H}^{n-1}.
\end{align*}
Then we get that
\begin{align*}
&\quad (-n-2\beta+2)r^{-n-2\beta+1}\int_{B_r(0)}(|\nabla
\vec{u}|^2-2F(\vec{u}))\\
&\quad+r^{-n-2\beta+2}\int_{\partial B_r(0)}(|\nabla
\vec{u}|^2-2F(\vec{u}))d\mathcal{H}^{n-1}\\
&\quad -\frac{\partial}{\partial r}\left(\beta
r^{-n-2\beta+1}\int_{\partial B_r(0)}
\vec{u}^2d\mathcal{H}^{n-1}\right)\\
&=2r^{-n-2\beta+2}\int_{\partial B_r(0)}(\nabla \vec{u}\cdot
\nu-\beta \frac{\vec{u}}{r})^2d\mathcal{H}^{n-1}\\
&\quad+2r^{-n-2\beta+1}\int_{B_r(0)}
\left(2(\beta-1)F(\vec{u})-\beta \vec{u}\vec{f}(\vec{u})\right),
\end{align*}
i.e.,
\begin{align}
(\Phi_{x_0}(r))'&=2r^{-n-2\beta+2}\int_{\partial B_r(0)}(\nabla
\vec{u}\cdot
\nu-\beta \frac{\vec{u}}{r})^2d\mathcal{H}^{n-1}\label{2.5}\\
&\quad+2r^{-n-2\beta+1}\int_{B_r(0)}
\left(2(\beta-1)F(\vec{u})-\beta \vec{u}\vec{f}(\vec{u})\right)\geq
0 \nonumber
\end{align}
 for a.e. $r\in (0,\delta)$. Integrating (\ref{2.5}) from $\rho$ to
$\sigma$, we can obtain (\ref{2.2}) and establish the monotonicity
formula in the theorem.
\end{proof}

Now we give some examples to illustrate Theorem \ref{emf05}.

{\bf Example 1.} Considering the following elliptic system (LES):
$$
\left\{ \begin{array}{ll} \Delta u+v=0 & {\rm in} \ \Omega,\\
\Delta v+u=0 & {\rm in} \ \Omega.
\end{array} \right.
$$
Then in this case, we know that $uf_1(u,v)=vf_2(u,v)=uv$ and
$F(u,v)=uv+c$, where $c$ is a real number. Assume that  $uv\in
L^1(B_{\delta}(x_0))$. Then for any $\beta$ such that
\begin{align}
\int_{B_r(x_0)}[(\beta-1)c-uv]\geq 0, \label{2.6}
\end{align}
we can get that
\begin{align*}
\Phi_{x_0}(r)&=r^{-n-2\beta+2}\int_{B_r(x_0)} (|\nabla u|^2+|\nabla
v|^2-2(uv+c))\\
&\quad-\beta r^{-n-2\beta+1}\int_{\partial B_r(x_0)}
(u^2+v^2)d\mathcal{H}^{n-1}
\end{align*}
is non-decreasing in $0<r<\delta$. In fact, we have that
\begin{align*}
(\Phi_{x_0}(r))'&=2r^{-n-2\beta+2}\int_{\partial B_r(x_0)}[(\nabla
u\cdot \nu-\beta \frac{u}{r})^2+(\nabla v\cdot
\nu-\beta \frac{v}{r})^2]d\mathcal{H}^{n-1}\\
&\quad +4r^{-n-2\beta+1}\int_{B_r(x_0)}[(\beta-1)c-uv]\geq 0
.\nonumber
\end{align*}
We often call (\ref{2.6}) the monotonicity condition for elliptic
system (LES). In particular, if $c>0$,  we can always take $\beta>1$
large enough; If $c<0$, we can take $\beta<1$ with $|\beta|$ large
enough; If $\beta=1$ and $uv\leq 0$ in
 $B_{\delta}(x_0)$, then the monotonicity condition is also true.

In particular, we restate our result when $u=v$. Assume that $u$ is
a variational solution of
$$
\Delta u+u=0 \quad {\rm in} \ \Omega
$$
and $B_{\delta}(x_0)\subset\subset \Omega\subset \subset
\mathbb{R}^n$. Then for $0<r<\delta$ and real constant $c$, we can
choose $\beta$ such that
$$
\Phi_{x_0}(r):=r^{-n-2\beta+2}\int_{B_r(x_0)} (|\nabla u|^2-u^2
-2c)-\beta r^{-n-2\beta+1}\int_{\partial B_r(x_0)}
u^2d\mathcal{H}^{n-1}
$$ is increasing in $r$ and satisfies
\begin{align*}
&\Phi_{x_0}(\sigma)-\Phi_{x_0}(\rho)=
\int_{\rho}^{\sigma}2r^{-n-2\beta+1}\int_{B_r(x_0)}(
2(\beta-1)c-u^2)\\
&\quad +\int_{\rho}^{\sigma}2r^{-n-2(\beta-1)}\int_{\partial
B_r(x_0)} (\nabla u\cdot \nu-\beta \frac{u}{r})^2
d\mathcal{H}^{n-1}dr\geq 0.
\end{align*}

{\bf Example 2.} We consider the famous Ginzburg-Landau model:
$$
\Delta \vec{u}+\frac{1}{\epsilon^2}\vec{u}(1-\vec{u}^2)=0, \quad
{\rm in} \ \Omega.
$$
Set
$$
F(\vec{u})=\frac{1}{4\epsilon^2}(1-\vec{u}^2)^2.
$$
Take $\beta>{1}$. Then
\begin{align*}
&2(\beta-1)F(\vec{u})-\beta\vec{u}f(\vec{u})\\
&=\frac{1}{2\epsilon^2}
[(\beta-1)-2(2\beta-1)\vec{u}^2+(3\beta-1)\vec{u}^4] \\
&\geq 0
\end{align*}
provided $\vec{u}^2 \leq \frac{\beta-1}{3\beta-1} $ or
$\vec{u}^2\geq 1$. Then
\begin{align*}
\Phi_{x_0}(r)&=r^{-n-2\beta+2}\int_{B_r(x_0)} (|\nabla \vec{u}|^2-
\frac{1}{4\epsilon^2}(1-\vec{u}^2)^2)\\
&\quad-\beta r^{-n-2\beta+1}\int_{\partial B_r(x_0)}
\vec{u}^2d\mathcal{H}^{n-1}.
\end{align*}

We have the following result

 \textbf{Assertion 2.1:} Let $\vec{u}\in
H^1_{loc}(\Omega)\cap L^4(\Omega)$ be a variational solution of
the Ginzburg-Landau model:
$$
\Delta \vec{u}+\frac{1}{\epsilon^2}\vec{u}(1-\vec{u}^2)=0, \quad
{\rm in} \ \Omega.
$$
Let $\beta>1$. Assume that
$$
\vec{u}^2 \leq \frac{\beta-1}{3\beta-1}\;\;{\rm or}\;\;
\vec{u}^2\geq 1
$$
in the ball $B_{\delta}(x_0)\subset \Omega$ for some $\delta>0$.
Then, we have for $0<r<\delta$ that
\begin{align*}
(\Phi_{x_0}(r))'&= r^{-n-2\beta+2}\int_{\partial B_r(x_0)}(\nabla
\vec{u}\cdot \nu-\beta \frac{u}{r})^2d\mathcal{H}^{n-1} \\
&+2r^{-n-2\beta+1}\int_{B_r(x_0)}\frac{1}{2\epsilon^2}
[(\beta-1)-2(2\beta-1)\vec{u}^2+(3\beta-1)\vec{u}^4]\\
&\geq 0.
\end{align*}

 \hfill $\Box$

{\bf Example 3.} Considering the elliptic system
$$
\left\{ \begin{array}{ll} \Delta u+\frac{u^pv^{q+1}}{q+1}=0 & {\rm in} \ \Omega,\\
\Delta v+\frac{u^{p+1}v^q}{p+1}=0 & {\rm in} \ \Omega,
\end{array} \right.
$$
where $(p,q)\geq 0$. Then, we know that
$$uf_1(u,v)=\frac{u^{p+1}v^{q+1}}{q+1}$$
and
 $$vf_2(u,v)=\frac{u^{p+1}v^{q+1}}{p+1}.$$
Take $$F(u,v)=\frac{u^{p+1}v^{q+1}}{(p+1)(q+1)}+c,$$ where $c$ is a
real number to be chosen. So,
\begin{align}
\Phi_{x_0}(r)&=r^{-n-2\beta+2}\int_{B_r(x_0)} (|\nabla u|^2+|\nabla
v|^2-2(\frac{u^{p+1}v^{q+1}}{(p+1)(q+1)}+c))\label{2.8}\\
&\quad-\beta r^{-n-2\beta+1}\int_{\partial B_r(x_0)}
(u^2+v^2)d\mathcal{H}^{n-1},\nonumber
\end{align}
and
\begin{align}
&(\Phi_{x_0}(r))'= r^{-n-2\beta+2}\int_{\partial B_r(0)}[(\nabla
u\cdot \nu-\beta \frac{u}{r})^2+(\nabla v\cdot
\nu-\beta \frac{v}{r})^2]d\mathcal{H}^{n-1}\label{2.9}\\
& +2r^{-n-2\beta+1}\int_{B_r(0)}\left(2(\beta-1)c-
\frac{\beta(p+q)+2}{(p+1)(q+1)}u^{p+1}v^{q+1} \right).\nonumber
\end{align}
As in example 1, we can choose suitable $c, \beta$ (see also Theorem
3.6 in (\cite{CPS}) for related stuff in parabolic case) such that
$$
\int_{B_r(x_0)}\left(2(\beta-1)c-
\frac{\beta(p+q)+2}{(p+1)(q+1)}u^{p+1}v^{q+1} \right)\geq 0
$$
and $\Phi_{x_0}(r)$ is increasing in $r$.
\hfill $\Box$

We now consider the blow-up (or blow-down) analysis for solutions to
(\ref{2.1}). Let $\vec{u}$ be a function in $B_{\delta}(x_0)$. For a
given point $x_0$ and a given sequence $\rho_k\rightarrow 0$, we
define the {\it scaled sequences} as follows
$$
\vec{u}_k(x):=\rho_k^{-\beta}\vec{u}(x_0+\rho_kx)
$$
and want to obtain more information on the solution's behavior. In
fact, we obtain the following theorem:

\begin{theorem}
Suppose that $0<\rho_k \rightarrow 0$ as $k\rightarrow \infty$, and
$\vec{u}$ is a solution of (\ref{2.1}) in $B_{\delta}(x_0)$ as in
Theorem 2.1, and that $\vec{u}$ satisfies at $x_0$ the growth
estimate
$$
{\rm sup}_{r\in (0,\delta)}(r^{-n-2\beta+1}\int_{\partial B_r(x_0)}
\vec{u}^2d\mathcal{H}^{n-1}+r^{-n-2\beta+2}|\int_{B_r(x_0)}2F(\vec{u})|)<\infty.
$$
Then as $r\rightarrow 0$, $\Phi(r)$ converges monotone
non-increasing to a limit, which is denoted by $M(u,x_0)$,  and for
any open $D\subset \subset \mathbb{R}^n$ and $k\geq k(D)$, the
scaled sequence $\vec{u}_k(x)$ is bounded in $H^{1}(D)$ and any weak
$H^{1}$-limit with respect to a subsequence $k\rightarrow \infty$ is
homogeneous of degree $\beta$.
\end{theorem}

{\bf Remark 2.1.} We say that the sequence
$\vec{u}_k(x):=\rho_k^{-\beta}\vec{u}(x_0+\rho_kx)$ is bounded in
$H^{1}(D)$ and any weak $H^{1}$-limit with respect to a subsequence
$k\rightarrow \infty$ is homogeneous of degree $\beta$, if every
$u_i$ is bounded in $H^{1}(D)$ and any weak $H^{1}$-limit with
respect to a subsequence $k\rightarrow \infty$ is homogeneous of
degree $\beta$.
\begin{proof}
First we can get for $0<R<\infty$ that
\begin{align*}
\Phi_{x_0}(\rho_k R)&=R^{-n-2\beta+2}\int_{B_R(0)}|\nabla
\vec{u}_k|^2-(\rho_k R)^{-n-2\beta+2}\int_{B_{\rho_kR}(0)}2F(\vec{u}))\\
&\quad-\beta R^{-n-2\beta+1}\int_{\partial
B_R(0)}\vec{u}^2_kd\mathcal{H}^{n-1},
\end{align*}
and we know that $\vec{u}_k$ is bounded in $H^{1}(D)$ for $k\geq
k(D)$ by the monotonicity formula and the condition of the
theorem.

By the results of Theorem 2.1, we know that $\Phi$ is non-decreasing
and bounded in $(0,r_0)$ for small and positive $r_0$, which means
that $\Phi$ has a right limit at $0$, and for $0<R<S$,
\begin{align*}
&0\leftarrow \Phi_{x_0}(\rho_kS)-\Phi_{x_0}(\rho_k R)\\
&\quad= \int^{\rho_k S }_{\rho_k
R}2r^{-n-2\beta+1}\int_{B_r(x_0)}(2(\beta-1)F(\vec{u})-\beta
\vec{u}\vec{f}(u))\\
&\qquad +\int^S_R 2r^{-n-2\beta+2}\int_{\partial B_r(0)} (\nabla
\vec{u_k}\cdot \nu-\beta \frac{\vec{u_k}}{r})^2
d\mathcal{H}^{n-1}dr\geq 0.
\end{align*}
i.e.
\begin{align*}
0&\leftarrow \Phi_{x_0}(\rho_kS)-\Phi_{x_0}(\rho_k R)\\
&\geq \int^S_R 2r^{-n-2\beta+2}\int_{\partial B_r(0)} (\nabla
\vec{u_k}\cdot \nu-\beta \frac{\vec{u_k}}{r})^2
d\mathcal{H}^{n-1}dr\\
&=\int_{B_S(0)\setminus B_R(0)} 2|x|^{-n-2\beta} (\nabla
\vec{u_k}(x)\cdot x-\beta \vec{u_k}(x))^2  \ {\rm as} \ k\rightarrow
\infty.
\end{align*}
Since the lower semi-continuity of the $L^2$-norm with respect to
weak convergence, we can take a subsequence $k\rightarrow \infty$
such that $\vec{u}_k\rightharpoonup \vec{u}_0$ weakly in
$H^{1}_{loc}(\mathbb{R}^n)$, then we obtain that $\nabla
\vec{u}_0(x)\cdot x-\beta \vec{u}_0(x)=0$ a.e. in $\mathbb{R}^n$.

From $\nabla \vec{u}_0(x)\cdot x-\beta \vec{u}_0(x)=0$, we can
easily prove that $\vec{u}_0$ is homogeneous of degree $\beta$.
\end{proof}

{\bf Remark 2.2.} If $\beta>0$, the scaled sequences are often
called {\it the blow-up sequences}. However, if $\beta<0$, our
scaled sequences are new.

\section{the monotonicity formula of a parabolic system}
In this section, we will consider the parabolic problem:
\begin{align}
\frac{\partial u_i}{\partial t}-\Delta u_i=f_i(u_1,...,u_m), \
i=1,...,m, \ {\rm in} \ (t_1,t_2)\times \mathbb{R}^n, \label{3.1}
\end{align}
where $t_1, t_2$ are two constants.

For convenience, we need some notations (see \cite{W3}).
Considering vector functions  $\vec{u}\in H^1_{loc}((0,T)\times
\mathbb{R}^n; \mathbb{R}^m)$ and $\psi\in H^1_{loc}((0,T)\times
\mathbb{R}^n;\mathbb{R}^{n+1})$, we denote by $\partial_t
\vec{u}:=\partial_0 \vec{u}$ the time derivative, by $\nabla
\vec{u}:=(\partial_1 \vec{u},..., \partial_n \vec{u})$ the space
gradient, by $\nabla_{t,x} \vec{u}:=(\partial_0 \vec{u},\partial_1
\vec{u},...,
\partial_n \vec{u})$ the time-space gradient, by ${\rm div}_{t,x}
\psi:=\sum^n_{k=0}
\partial_k \psi_k$ the time-space divergence, and by
$$
D\psi:=\left(\begin{array}{ll} \partial_1 \psi_0 \qquad  ... \
&\partial_n \psi_0\\
 \partial_1 \psi_1 \qquad  ... \
&\partial_n \psi_1\\
\qquad \qquad ...& \\
\partial_1 \psi_n \quad ... \ &\partial_n \psi_n
\end{array} \right)
$$
the space Jacobian. Moreover, denote by $(\nabla \vec{u}\cdot
x+2t\partial_t \vec{u}-\beta \vec{u})^2
 =\sum_{i=1}^m (\nabla u_i\cdot x+2t\partial_t u_i-\beta u_i)^2$.

Next, we give the definition of a variational solution of
(\ref{3.1}).

{\bf Definition 2.} {\it We define $\vec{u}\in H^1((t_1,t_2)\times
B_R(0))$ for any $R\in (0,\infty)$ to be a variational solution of
(\ref{3.1}) if $\vec{u}\in H^1((t_1,t_2)\times \mathbb{R}^n)$
satisfies (\ref{3.1}) in the distributional sense with
\begin{equation}
u_if_i(\vec{u}),\ F(\vec{u})\in L^1((t_1,t_2)\times
\mathbb{R}^n) \label{3.c}
\end{equation}
for $i=1,...,m$ and the vanishing of first variation at
$\vec{v}=\vec{u}$ with respect to variations of the domain in time
and spaces of the following functional
$$
\mathbb{G}(\vec{u},\vec{v}):=\int^{t_2-\delta}_{t_1+\delta}
\int_{\mathbb{R}^n}(|\nabla \vec{v}|^2-2
F(\vec{v}))+\int^{t_2-\delta}_{t_1+\delta}\int_{\mathbb{R}^n}2\vec{v}
\partial_t \vec{u},
$$
where $(\vec{u},\vec{v})\in (H^1)^2$ and $\vec{v}$ also satisfies
(\ref{3.c}), i.e.
$$\frac{d}{d\varepsilon}\mathbb{G}(\vec{u},\vec{u}((t,x)+\varepsilon
\psi(t,x)))|_{\varepsilon=0}=0,$$ that is,
\begin{align*}
0&=\int^{t_2-\delta}_{t_1+\delta}\int_{\mathbb{R}^n}[(|\nabla
\vec{u}|^2-2F(\vec{u})){\rm div}_{t,x}\psi-2\nabla_{t,x}
\vec{u}D\psi \nabla
\vec{u}-2\partial_t \vec{u} \nabla_{t,x}\vec{u}\cdot \psi]\\
&\quad -[\int_{\mathbb{R}^n}(|\nabla
\vec{u}|^2-2F(\vec{u}))\psi_0]^{t_2-\delta}_{t_1+\delta}\\
&=\int^{t_2-\delta}_{t_1+\delta}\int_{\mathbb{R}^n}[(|\nabla
\vec{u}|^2-2F(\vec{u}))\sum_{k=0}^n\partial_k
\psi_k-2\sum_{j=1}^n\sum_{k=0}^n \partial_j
\vec{u}\partial_j\psi_k\partial_k \vec{u}\\
&\quad-2\partial_t \vec{u}\sum_{k=0}^n\partial_k
\vec{u}\psi_k]-\int_{\mathbb{R}^n}[(|\nabla
\vec{u}|^2-2F(\vec{u}))\psi_0](t_2-\delta)\\
&\quad+\int_{\mathbb{R}^n}[(|\nabla
\vec{u}|^2-2F(\vec{u}))\psi_0](t_1+\delta)
\end{align*}
for a.e. small and positive $\delta$ and any $\psi\in
C^1(\mathbb{R}^{n+1})$ such that
$$supp\psi(t)\subset\subset
\mathbb{R}^n$$ for any $t\in(t_1,t_2)$. }

We now state a monotonicity formula for variational solution of
(\ref{3.1}).

\begin{theorem}
(monotonicity formula). Let  $\vec{u}$ be a variational solution
of (\ref{3.1}) in $((t_1,T)\cup (T,t_2))\times \mathbb{R}^n$,
where $t_1\leq T\leq t_2$. Let $x_0\in \mathbb{R}^n$. Assume that
\begin{align*}
&\quad {\rm sup}_{t\in (t_1,
T-\delta)\cup(T+\delta,t_2)}\int_{\mathbb{R}^n}\exp(-\frac{|x-x_0|^2}{4(T-t)})(|\nabla
\vec{u}|^2-2F(\vec{u}))(t,x)dx\\
&+\int_{(t_1,T-\delta)\cup(T+\delta,t_2)}\int_{\mathbb{R}^n}\exp(-\frac{|x-x_0|^2}{4(T-t)})
((\partial_t \vec{u})^2+\vec{u}^2)(t,x)dxdt<\infty
\end{align*}
for any small positive $\delta$. Then for any real $\beta$ such that
$$\int_{T^-_r} G_{(T,x_0)}[2(\beta-1)F(\vec{u})-\beta \vec{u}\vec{f}(\vec{u})]\geq 0$$
and
$$\int_{T^+_r} G_{(T,x_0)}[2(\beta-1)F(\vec{u})-\beta \vec{u}\vec{f}(\vec{u})]\geq 0,$$
the functions
\begin{align*}
\Psi^-(r)&:=r^{-2\beta}\int_{T^-_r(T)}(|\nabla
\vec{u}|^2-2F(\vec{u}))G_{(T,x_0)}\\
&\quad-\frac{\beta}{2}r^{-2\beta}\int_{T^-_r(T)}\frac{1}{T-t}\vec{u}^2G_{(T,x_0)}
\end{align*}
and
\begin{align*}
\Psi^+(r)&:=r^{-2\beta}\int_{T^+_r(T)}(|\nabla
\vec{u}|^2-2F(\vec{u}))G_{(T,x_0)}\\
&\quad-\frac{\beta}{2}r^{-2\beta}\int_{T^+_r(T)}\frac{1}{T-t}\vec{u}^2G_{(T,x_0)}
\end{align*}
are well defined in the interval $(0,\frac{\sqrt{T-t_1}}{2})$ and
$(0,\frac{\sqrt{t_2-T}}{2})$, respectively, and they satisfy for any
$0<\rho<\sigma<\frac{\sqrt{T-t_1}}{2}$ and
$0<\rho<\sigma<\frac{\sqrt{t_2-T}}{2}$, respectively, the
monotonicity formulae
\begin{align}
&\Psi^-(\sigma)-\Psi^-(\rho)=\int^{\sigma}_{\rho}2r^{-2\beta-1}
\int_{T^-_r(T)}[2(\beta-1)F(\vec{u})-\beta
\vec{u}\vec{f}(\vec{u})]G_{(T,x_0)} \label{3.2} \\
&+\int^{\sigma}_{\rho}r^{-2\beta-1}\int_{T^-_r(T)}\frac{1}{T-t}(\nabla
\vec{u}\cdot(x-x_0)-2(T-t)\partial_t \vec{u}-\beta \vec{u})^2G_{(T,x_0)}\nonumber\\
&\geq 0\nonumber
\end{align}
and
\begin{align}
&\Psi^+(\sigma)-\Psi^+(\rho)=\int^{\sigma}_{\rho}2r^{-2\beta-1}
\int_{T^+_r(T)}[2(\beta-1)F(\vec{u})-\beta
\vec{u}\vec{f}(\vec{u})]G_{(T,x_0)}\label{3.3}\\
&+\int^{\sigma}_{\rho}r^{-2\beta-1}\int_{T^+_r(T)}\frac{1}{T-t}(\nabla
\vec{u}\cdot(x-x_0)-2(T-t)\partial_t \vec{u}-\beta \vec{u})^2G_{(T,x_0)}\nonumber\\
&\geq 0\nonumber.
\end{align}
\end{theorem}

\begin{proof}
We will only give a proof for the monotonicity of $\Psi^-$ because
we can replace in what follows the interval $(-4r^2,-r^2)$ by
$(r^2,4r^2)$ in order to obtain a proof for $\Psi^+$. Without loss
of generality, we can assume that $x_0=0$ and $T=0$. We omit the
index $(0,0)$ in $G_{(0,0)}$ and simply denote it by $G$, and
denote $T^-_r(0)$ by $T^-_r$. Choosing $t_1:=-4r^2$, $t_2:=-r^2$,
and $\psi(t,x):=(2t,x)G(t,x)\eta_{\epsilon}(x)$ in Definition 2
where $\eta_{\epsilon}\in H^{1,\infty}(\mathbb{R}^n)$ will be
chosen later, we obtain that
\begin{align}
0&=\int_{T^-_r}[(|\nabla
\vec{u}|^2-2F(\vec{u}))(2G+2t\partial_tG+{\rm
div}(xG))\eta_{\epsilon}\label{3.4}\\
&\quad-2\eta_{\epsilon}\sum^n_{j=1}\sum^n_{k=1}\partial_j\vec{u}
(\delta_{jk}G+\partial_jGx_k)\partial_k
\vec{u}-2\eta_{\epsilon}\sum^n_{j=1}\partial_j
\vec{u}\partial_j G2t\partial_t\vec{u}\nonumber\\
&\qquad-2\eta_{\epsilon}\sum^n_{j=1}\partial_j
\vec{u}Gx_j\partial_t
\vec{u}-2\eta_{\epsilon}(\partial_t \vec{u})^22tG]\nonumber\\
&\quad-\int_{\mathbb{R}^n}[2t\eta_{\epsilon}(|\nabla
\vec{u}|^2-2F(\vec{u}))G](-r^2)\nonumber\\
&\quad+\int_{\mathbb{R}^n}[2t\eta_{\epsilon}(|\nabla
\vec{u}|^2-2F(\vec{u}))G](-4r^2)\nonumber\\
&\quad+\int_{T^-_r}[(|\nabla \vec{u}|^2-2F(\vec{u}))\nabla \eta_{\epsilon}\cdot xG\nonumber\\
&\quad-2\sum^n_{j=1}\sum^n_{k=1}\partial_j \vec{u}\partial_j
\eta_{\epsilon}x_k\partial_k
\vec{u}G-2\sum^n_{j=1}\partial_j\vec{u}\partial_j\eta_{\epsilon}2t\partial_t\vec{u}G]\nonumber
\end{align}
for a.e. $r\in (0, \frac{\sqrt{T-t_1}}{2})$.

Multiplying (\ref{3.4}) by $-r^{-2\beta-1}$ and choosing
$\eta_{\epsilon}(x):=\min(1,\max(0,2-\epsilon|x|))$ for small
$\epsilon>0$, we get that
\begin{align*}
0&=r^{-2\beta-1}[\int_{\mathbb{R}^n}2t\eta_{\epsilon}G(|\nabla
\vec{u}|^2-2F(\vec{u}))]^{(-r^2)}_{(-4r^2)}\\
&\quad-2\beta r^{-2\beta-1}\int_{T^-_r}\eta_{\epsilon}G(\nabla
\vec{u}|^2-2F(\vec{u}))\nonumber\\
&\quad+2\beta r^{-2\beta-1}\int_{T^-_r}\eta_{\epsilon}G|\nabla
\vec{u} |^2-4(\beta-1)r^{-2\beta-1}
\int_{T^-_r}\eta_{\epsilon}GF(\vec{u})\nonumber\\
&\quad+r^{-2\beta-1}\int_{T^-_r} \frac{\eta_{\epsilon}G}{t} (\nabla
\vec{u}\cdot x)^2+
r^{-2\beta-1}\int_{T^-_r} 4\eta_{\epsilon}G\nabla \vec{u}\cdot x\partial_t \vec{u}\nonumber\\
&\quad+r^{-2\beta-1}\int_{T^-_r}[4\eta_{\epsilon}Gt(\partial_t\vec{u})^2-
(|\nabla \vec{u}|^2-2F(\vec{u}))\nabla \eta_{\epsilon}\cdot xG]\nonumber\\
&\quad+r^{-2\beta-1}\int_{T^-_r}2[G\nabla \vec{u} \cdot \nabla
\eta_{\epsilon} \nabla \vec{u}\cdot x+\nabla \vec{u}\cdot \nabla
\eta_{\epsilon} 2t
\partial_t \vec{u} G],\nonumber
\end{align*}
where we use the fact that $\nabla G=\frac{xG}{2t}$ and $\partial_t
G+\Delta G=0$ in $\{t<0\}\cup \{t>0\}$.

As in the proof of (\ref{2.4}), we obtain that
\begin{align}
\int_{T^-_r}|\nabla \vec{u}|^2G\eta_{\epsilon}&=-\int_{T^-_r}
[\vec{u}\eta_{\epsilon}\nabla \vec{u}\cdot \nabla
G+\eta_{\epsilon}G\vec{u}(\partial_t\vec{u}-\vec{f}(\vec{u}))\label{3.5}\\
&\qquad \qquad \qquad +\vec{u}G\nabla\eta_{\epsilon}\cdot\nabla
\vec{u}]. \nonumber
\end{align}
Using (\ref{3.5}), we can get that
\begin{align*}
0&=r^{-2\beta-1}[\int_{\mathbb{R}^n}2t\eta_{\epsilon}G(|\nabla
\vec{u}|^2-2F(\vec{u}))]^{(-r^2)}_{(-4r^2)}\\
&\quad-2\beta r^{-2\beta-1}\int_{T^-_r}\eta_{\epsilon}G(\nabla
\vec{u}|^2-2F(\vec{u}))\nonumber\\
&\quad+2\beta
r^{-2\beta-1}\int_{T^-_r}[\eta_{\epsilon}G\vec{u}(\vec{f}(\vec{u})-\partial_t
\vec{u})-\frac{\eta_{\epsilon}G\vec{u}}{2t} \nabla \vec{u}\cdot x]\nonumber\\
&\quad-2\beta r^{-2\beta-1}\int_{T^-_r}\vec{u}G\nabla
\eta_{\epsilon}\cdot
\nabla \vec{u}-4(\beta-1)r^{-2\beta-1}\int_{T^-_r}\eta_{\epsilon}GF(\vec{u})\nonumber\\
&\quad+r^{-2\beta-1}\int_{T^-_r}[\frac{\eta_{\epsilon}G}{t}(\nabla
\vec{u}\cdot x+2t\partial_t \vec{u})^2-
(|\nabla \vec{u}|^2-2F(\vec{u}))\nabla \eta_{\epsilon}\cdot xG]\nonumber\\
&\quad+r^{-2\beta-1}\int_{T^-_r}2[G\nabla \vec{u} \cdot \nabla
\eta_{\epsilon}
\nabla \vec{u}\cdot x+\nabla \vec{u}\cdot \nabla \eta_{\epsilon} 2t \partial_t \vec{u} G]\nonumber\\
\end{align*}
i.e.
\begin{align}
0&=r^{-2\beta-1}[\int_{\mathbb{R}^n}2t\eta_{\epsilon}G(|\nabla
\vec{u}|^2-2F(\vec{u}))]^{(-r^2)}_{(-4r^2)} \label{3.6} \\
&\quad-2\beta r^{-2\beta-1}\int_{T^-_r}\eta_{\epsilon}G(\nabla
\vec{u}|^2-2F(\vec{u}))\nonumber\\
&\quad+r^{-2\beta-1}\int_{T^-_r}[\frac{\eta_{\epsilon}G}{t}(\nabla
\vec{u}\cdot
x+2t\partial_t \vec{u}-\beta \vec{u})^2\nonumber\\
&\quad +2r^{-2\beta-1}\int_{T^-_r}(\beta
\vec{u}\vec{f}(\vec{u})-2(\beta-1)F(\vec{u}))\eta_{\epsilon}G\nonumber\\
&\quad+r^{-2\beta-1}\int_{T^-_r}\eta_{\epsilon}G(\frac{\beta}{t}\vec{u}\nabla
\vec{u}\cdot x+2\beta \vec{u}\partial_t
\vec{u}-\frac{\beta^2}{t}\vec{u}^2)\nonumber\\
&\quad +o(1)  \nonumber
\end{align}
as $ \epsilon\rightarrow 0$. Notice that
\begin{align*}
&\quad
r^{-2\beta-1}\int_{T^-_r}\eta_{\epsilon}G(\frac{\beta}{t}\vec{u}\nabla
\vec{u}\cdot
x+2\beta \vec{u}\partial_t \vec{u}-\frac{\beta^2}{t}\vec{u}^2)\\
&=\beta
\int_{T^-_1}\eta_{\epsilon}(rx)\frac{G(t,x)}{t}r^{-\beta}u(r^2t,rx)
[r^{-\beta}(\nabla \vec{u})(r^2t,rx)\cdot x\\
&\quad +r^{-\beta+1}2t(\partial_t\vec{u})(r^2t,rx)-\beta
r^{-\beta-1}\vec{u}(r^2t,rx)]\\
&=o(1)+\partial_r(\frac{\beta}{2}\int_{T^-_1}(\frac{\vec{u}(r^2t,rx)}
{r^{\beta}})^2\frac{G(t,x)}{t}\eta_{\epsilon}(rx))\\
&=o(1)+\partial_r(\frac{\beta}{2}\int_{T^-_r}(\frac{\vec{u}(x,t)}
{r^{\beta}})^2\frac{G(t,x)}{t}\eta_{\epsilon}(x))
\end{align*}
as $ \epsilon\rightarrow 0$. Letting $\epsilon\rightarrow 0$ in
(\ref{3.6}), we find that
\begin{align}
&(\Psi^-(r))'=2r^{-2\beta-1}
\int_{T^-_r(T)}[2(\beta-1)F(\vec{u})-\beta
\vec{u}\vec{f}(\vec{u})]G_{(T,x_0)}\label{3.7}\\
&+r^{-2\beta-1}\int_{T^-_r(T)}\frac{1}{T-t}(\nabla
\vec{u}\cdot(x-x_0)-2(T-t)\partial_t \vec{u}-\beta \vec{u})^2G_{(T,x_0)}\nonumber\\
&\geq 0\nonumber.
\end{align}
Integrating (\ref{3.7}) from $\rho$ to $\sigma$, we can obtain
(\ref{3.2}). The proof of the theorem is complete.
\end{proof}

Now let's consider an example. Assume that $t_1<T<t_2$, $x_0\in
\mathbb{R}^n$, and
$$u\in H^1( ((t_1,T)\cup (T,t_2))\times \mathbb{R}^n)$$ is a
variational solution of
$$u_t-\Delta u=u \quad {\rm in} \ (t_1, t_2)\times
\mathbb{R}^n.$$ Suppose furthermore that
\begin{align*}
&\quad {\rm sup}_{t\in (t_1,
T-\delta)\cup(T+\delta,t_2)}\int_{\mathbb{R}^n}\exp(-\frac{|x-x_0|^2}{4(T-t)})(|\nabla
u|^2-u^2)dx\\
&+\int_{(t_1,T-\delta)\cup(T+\delta,t_2)}\int_{\mathbb{R}^n}\exp(-\frac{|x-x_0|^2}{4(T-t)})
((\partial_t u)^2+u^2)(t,x)dxdt<\infty
\end{align*}
for any positive $\delta$. Then there exist two constants $\beta,c$
such that the functions
\begin{align*}
\Psi^-(r)&:=r^{-2\beta}\int_{T^-_r(T)}(|\nabla
u|^2-u^2-2c)G_{(T,x_0)}\\
&\quad-\frac{\beta}{2}r^{-2\beta}\int_{T^-_r(T)}\frac{1}{T-t}u^2G_{(T,x_0)}
\end{align*}
and
\begin{align*}
\Psi^+(r)&:=r^{-2\beta}\int_{T^+_r(T)}(|\nabla
u|^2-u^2-2c)G_{(T,x_0)}\\
&\quad-\frac{\beta}{2}r^{-2\beta}\int_{T^+_r(T)}\frac{1}{T-t}u^2G_{(T,x_0)}
\end{align*}
are well defined in the interval $(0,\frac{\sqrt{T-t_1}}{2})$ and
$(0,\frac{\sqrt{t_2-T}}{2})$, respectively, and satisfy for any
$0<\rho<\sigma<\frac{\sqrt{T-t_1}}{2}$ and
$0<\rho<\sigma<\frac{\sqrt{t_2-T}}{2}$, respectively, the
monotonicity formulae
\begin{align*}
&\Psi^-(\sigma)-\Psi^-(\rho)=\int^{\sigma}_{\rho}2r^{-2\beta-1}
\int_{T^-_r(T)}(2(\beta-1)c-u^2)\\
&+\int^{\sigma}_{\rho}r^{-2\beta-1}\int_{T^-_r(T)}\frac{1}{T-t}(\nabla
u\cdot(x-x_0)-2(T-t)\partial_t u-\beta u)^2G_{(T,x_0)}\\
&\geq 0
\end{align*}
and
\begin{align*}
&\Psi^+(\sigma)-\Psi^+(\rho)=\int^{\sigma}_{\rho}2r^{-2\beta-1}
\int_{T^+_r(T)}(2(\beta-1)c-u^2)\\
&+\int^{\sigma}_{\rho}r^{-2\beta-1}\int_{T^+_r(T)}\frac{1}{T-t}(\nabla
u\cdot(x-x_0)-2(T-t)\partial_t u-\beta u)^2G_{(T,x_0)}\\
&\geq 0.
\end{align*}

In the remaining part of this section, we will consider the free
boundary problem:
$$
(\star) \left\{ \begin{array}{ll} \frac{\partial u_i}{\partial
t}-\Delta u_i=\chi_{\Omega} f_i(\vec{u}) & {\rm
in} \ \mathbb{R}^n\times \mathbb{R},\\
u_i=|\nabla u_i|=0, i=1,...,m, & {\rm in} \ (\mathbb{R}^n\times
\mathbb{R}) \setminus \Omega,
\end{array} \right.
$$
where $$ \Lambda:=\{u_i=|\nabla u_i|=0\}, \quad \Omega:=
(\mathbb{R}^n\times \mathbb{R})\setminus \Lambda. $$

As before, we start with a definition of variational solution to the
problem above.

 {\bf Definition 3.} {\it We define
$\vec{u}\in H^1((t_1,t_2)\times B_R(0))$ for any $R\in (0,\infty)$
to be a variational solution of $(\star)$, if $\vec{u}\in
H^1((t_1,t_2)\times \mathbb{R}^n)$ satisfies ($\star$) in the
distributional sense with
$$u_if_i(\vec{u}),\ F(\vec{u})\in L^1((t_1,t_2)\times \mathbb{R}^n)
$$
for $i=1,...,m$ and the first variation with respect to variations
of the domain in time and spaces of the functional
$$
\mathbb{G}(\vec{u},\vec{v}):=\int^{t_2-\delta}_{t_1+\delta}
\int_{\mathbb{R}^n}(|\nabla \vec{v}|^2-2\chi_{\Omega}F(\vec{v}))+
\int^{t_2-\delta}_{t_1+\delta}
\int_{\mathbb{R}^n}2\vec{v}\partial_t \vec{u},
$$
vanishes at $\vec{v}=\vec{u}$, i.e.
$\frac{d}{d\varepsilon}\mathbb{G}(\vec{u},\vec{u}((t,x)+\varepsilon
\psi(t,x)))|_{\varepsilon=0}=0$, that is
\begin{align*}
0&=\int^{t_2-\delta}_{t_1+\delta}\int_{\mathbb{R}^n}[(|\nabla
\vec{u}|^2-2\chi_{\Omega}F(\vec{u})){\rm
div}_{t,x}\psi-2\nabla_{t,x} \vec{u}D\psi \nabla
\vec{u}-2\partial_t \vec{u} \nabla_{t,x}\vec{u}\cdot \psi]\\
&\quad -[\int_{\mathbb{R}^n}(|\nabla
\vec{u}|^2-2\chi_{\Omega}F(\vec{u}))\psi_0]^{t_2-\delta}_{t_1+\delta}\\
&=\int^{t_2-\delta}_{t_1+\delta}\int_{\mathbb{R}^n}[(|\nabla
\vec{u}|^2-2\chi_{\Omega}F(\vec{u}))\sum_{k=0}^n\partial_k
\psi_k-2\sum_{j=1}^n\sum_{k=0}^n \partial_j
\vec{u}\partial_j\psi_k\partial_k \vec{u}\\
&\quad-2\partial_t \vec{u}\sum_{k=0}^n\partial_k
\vec{u}\psi_k]-\int_{\mathbb{R}^n}[(|\nabla
\vec{u}|^2-2\chi_{\Omega}F(\vec{u}))\psi_0](t_2-\delta)\\
&\quad+\int_{\mathbb{R}^n}[(|\nabla
\vec{u}|^2-2\chi_{\Omega}F(\vec{u}))\psi_0](t_1+\delta)
\end{align*}
for a.e. small and positive $\delta$ and any $\psi\in C^1(
\mathbb{R}^n\times \mathbb{R})$ such that {\rm
supp}$\psi(t)\subset\subset \mathbb{R}^n$ for any $t\in(t_1,t_2)$.
}

Then we have the following result.
\begin{theorem}
Assume that $\vec{u}$ is a variational solution of $(\star)$,
$(x_0,T)\in \Lambda$, $t_1\leq T\leq t_2$, and
\begin{align*}
&\quad {\rm sup}_{t\in (t_1,
T-\delta)\cup(T+\delta,t_2)}\int_{\mathbb{R}^n}\exp(-\frac{|x-x_0|^2}{4(T-t)})(|\nabla
\vec{u}|^2-2\chi_{\Omega}F(\vec{u}))(t,x)dx\\
&+\int_{(t_1,T-\delta)\cup(T+\delta,t_2)}\int_{\mathbb{R}^n}\exp(-\frac{|x-x_0|^2}{4(T-t)})
((\partial_t \vec{u})^2+\vec{u}^2)(t,x)dxdt<\infty
\end{align*}
for any positive $\delta$. Let $\varphi(x)\geq 0$ be a
$C^{\infty}$ cut-off function in $\mathbb{R}^n$ with ${\rm
supp}\varphi\subset B_{3/4}(x_0)$ and $\varphi|_{B_{1/2}(x_0)}=1$.

Then for any $\beta$ such that
$$\int_{T^-_r} \chi_{\Omega}[2(\beta-1)F(\vec{u})-
\beta \vec{u}\vec{f}(\vec{u})]\varphi G_{(T,x_0)}\geq 0$$ and
$$\int_{T^+_r} \chi_{\Omega}[2(\beta-1)F(\vec{u})-
\beta \vec{u}\vec{f}(\vec{u})]\varphi G_{(T,x_0)}\geq 0,$$ there
exists constant $C=C(n, \varphi, \beta)>0$ such that  the functions
\begin{align*}
\Psi^-(r)&:=r^{-2\beta}\int_{T^-_r(T)}(|\nabla
\vec{u}|^2-2\chi_{\Omega}F(\vec{u}))\varphi G_{(T,x_0)}\\
&\quad-\frac{\beta}{2}r^{-2\beta}\int_{T^-_r(T)}\frac{1}{T-t}\vec{u}^2\varphi
G_{(T,x_0)}+C\int^r_0 s^{-n-2\beta-1}e^{-1/(16s^2)}ds
\end{align*}
and
\begin{align*}
\Psi^+(r)&:=r^{-2\beta}\int_{T^+_r(T)}(|\nabla
\vec{u}|^2-2\chi_{\Omega}F(\vec{u})) \varphi G_{(T,x_0)}\\
&\quad-\frac{\beta}{2}r^{-2\beta}\int_{T^+_r(T)}\frac{1}{T-t}\vec{u}^2\varphi
G_{(T,x_0)}+C\int^r_0 s^{-n-2\beta-1}e^{-1/(16s^2)}ds
\end{align*}
are well defined in the interval $(0,\frac{\sqrt{T-t_1}}{2})$ and
$(0,\frac{\sqrt{t_2-T}}{2})$, respectively, and satisfy for any
$0<\rho<\sigma<\frac{\sqrt{T-t_1}}{2}$ and
$0<\rho<\sigma<\frac{\sqrt{t_2-T}}{2}$, respectively, the
monotonicity formulae
\begin{align*}
&\Psi^-(\sigma)-\Psi^-(\rho)=\int^{\sigma}_{\rho}2r^{-2\beta-1}
\int_{T^-_r(T)}\chi_{\Omega}[2(\beta-1)F(\vec{u})-\beta
\vec{u}\vec{f}(\vec{u})]
\varphi G_{(T,x_0)}\\
&+\int^{\sigma}_{\rho}r^{-2\beta-1}\int_{T^-_r(T)}\frac{1}{T-t}(\nabla
\vec{u}\cdot(x-x_0)-2(T-t)\partial_t \vec{u}-\beta \vec{u})^2\varphi G_{(T,x_0)}\\
&+C\int^{\sigma}_{\rho} s^{-n-2\beta-1}e^{-1/(16s^2)}ds\geq 0
\end{align*}
and
\begin{align*}
&\Psi^+(\sigma)-\Psi^+(\rho)=\int^{\sigma}_{\rho}2r^{-2\beta-1}
\int_{T^+_r(T)}\chi_{\Omega}[2(\beta-1)F(\vec{u})-\beta
\vec{u}\vec{f}(\vec{u})]
\varphi G_{(T,x_0)}\\
&+\int^{\sigma}_{\rho}r^{-2\beta-1}\int_{T^+_r(T)}\frac{1}{T-t}(\nabla
\vec{u}\cdot(x-x_0)-2(T-t)\partial_t \vec{u}-\beta \vec{u})^2\varphi G_{(T,x_0)}\\
&+C\int^{\sigma}_{\rho} s^{-n-2\beta-1}e^{-1/(16s^2)}ds\geq 0.
\end{align*}
\end{theorem}

We remark that the proof is similar to that of Theorem 3.1. However,
for completeness, we give a full proof.

\begin{proof}
As before, we only give a proof for the monotonicity of $\Psi^-$,
because we can replace in what follows the interval $(-4r^2,-r^2)$
by $(r^2,4r^2)$ in order to obtain a proof with respect to
$\Psi^+$. Without loss of generality, we can assume that $x_0=0$
and $T=0$. We omit the index $(0,0)$ of $G_{(0,0)}$ and denote it
by $G$, and denote $T^-_r(0,0)$ by $T^-_r$. Choosing $t_1:=-4r^2$,
$t_2:=-r^2$, and $\psi(t,x):=(2t,x)G(t,x)\varphi(x)$ in Definition
3, where $\varphi(x)$ is a $C^{\infty}$ cut-off function in
$\mathbb{R}^n$ with ${\rm supp}\varphi\subset B_{3/4}(x_0)$ and
$\varphi|_{B_{1/2}(x_0)}=1$, we obtain that
\begin{align}
0&=\int_{T^-_r}[(|\nabla
\vec{u}|^2-2\chi_{\Omega}F(\vec{u}))(2G+2t\partial_tG+{\rm
div}(xG))\varphi\label{3.8}\\
&\quad-2\varphi\sum^n_{j=1}\sum^n_{k=1}\partial_j\vec{u}
(\delta_{jk}G+\partial_jGx_k)\partial_t\vec{u}-2\varphi\sum^n_{j=1}\partial_j
\vec{u}\partial_j G2t\partial_t\vec{u}\nonumber\\
&\qquad-2\varphi\sum^n_{j=1}\partial_j \vec{u}Gx_j\partial_t
\vec{u}-2\varphi(\partial_t \vec{u})^22tG]\nonumber\\
&\quad-\int_{\mathbb{R}^n}[2t\varphi(|\nabla
\vec{u}|^2-2\chi_{\Omega}F(\vec{u}))G](-r^2)\nonumber\\
&\quad+\int_{\mathbb{R}^n}[2t\varphi(|\nabla
\vec{u}|^2-2\chi_{\Omega}F(\vec{u}))G](-4r^2)\nonumber\\
&\quad+\int_{T^-_r}[(|\nabla \vec{u}|^2-2\chi_{\Omega}F(\vec{u}))\nabla \varphi\cdot xG\nonumber\\
&\quad-2\sum^n_{j=1}\sum^n_{k=1}\partial_j \vec{u}\partial_j
\varphi x_k\partial_k
\vec{u}G-2\sum^n_{j=1}\partial_j\vec{u}\partial_j\varphi2t\partial_t\vec{u}G]\nonumber
\end{align}
for a.e. $r\in (0, \frac{\sqrt{T-t_1}}{2})$.

Multiplying (\ref{3.8}) by $-r^{-2\beta-1}$, we get that
\begin{align*}
0&=r^{-2\beta-1}[\int_{\mathbb{R}^n}2t\varphi G(|\nabla
\vec{u}|^2-2\chi_{\Omega}F(\vec{u}))]^{(-r^2)}_{(-4r^2)}\\
&\quad-2\beta r^{-2\beta-1}\int_{T^-_r}\varphi G(|\nabla
\vec{u}|^2-2\chi_{\Omega}F(\vec{u}))\nonumber\\
&\quad+2\beta r^{-2\beta-1}\int_{T^-_r}\varphi G|\nabla \vec{u}
|^2-4(\beta-1)r^{-2\beta-1}
\int_{T^-_r}\chi_{\Omega}\varphi GF(\vec{u})\nonumber\\
&\quad+r^{-2\beta-1}\int_{T^-_r} \frac{\varphi G}{t} (\nabla
\vec{u}\cdot x)^2+
r^{-2\beta-1}\int_{T^-_r} 4\varphi G\nabla \vec{u}\cdot x\partial_t \vec{u}\nonumber\\
&\quad+r^{-2\beta-1}\int_{T^-_r}[4\varphi Gt(\partial_t\vec{u})^2-
(|\nabla \vec{u}|^2-2\chi_{\Omega}F(\vec{u}))\nabla \varphi\cdot xG]\nonumber\\
&\quad+r^{-2\beta-1}\int_{T^-_r}2[G\nabla \vec{u} \cdot \nabla
\varphi \nabla \vec{u}\cdot x+\nabla \vec{u}\cdot \nabla \varphi 2t
\partial_t \vec{u} G],
\end{align*}
where we use the fact that $\nabla G=\frac{xG}{2t}$ and $\partial_t
G+\Delta G=0$ in $\{t<0\}\cup \{t>0\}$.

As in the proof of (\ref{3.6}), we obtain that
\begin{align}
\int_{T^-_r}|\nabla \vec{u}|^2G\varphi&=-\int_{T^-_r}
[\vec{u}\varphi\nabla \vec{u}\cdot \nabla
G+\varphi G\vec{u}(\partial_t\vec{u}-\chi_{\Omega}\vec{f}(\vec{u}))\nonumber\\
&\qquad \qquad \qquad +\vec{u}G\nabla\varphi\cdot\nabla
\vec{u}].\label{3.9}
\end{align}
Using (\ref{3.9}), we can get that
\begin{align*}
0&=r^{-2\beta-1}[\int_{\mathbb{R}^n}2t\varphi G(|\nabla
\vec{u}|^2-2\chi_{\Omega}F(\vec{u}))]^{(-r^2)}_{(-4r^2)}\nonumber\\
&\quad-2\beta r^{-2\beta-1}\int_{T^-_r}\varphi G(|\nabla
\vec{u}|^2-2\chi_{\Omega}F(\vec{u}))\nonumber\\
&\quad+2\beta r^{-2\beta-1}\int_{T^-_r}[\varphi
G\vec{u}(\chi_{\Omega}\vec{f}(\vec{u})-\partial_t
\vec{u})-\frac{\varphi G\vec{u}}{2t} \nabla \vec{u}\cdot x]\nonumber\\
&\quad-2\beta r^{-2\beta-1}\int_{T^-_r}\vec{u}G\nabla \varphi \cdot
\nabla \vec{u}-4(\beta-1)r^{-2\beta-1}\int_{T^-_r}\chi_{\Omega}\varphi GF(\vec{u})\nonumber\\
&\quad+r^{-2\beta-1}\int_{T^-_r}[\frac{\varphi G}{t}(\nabla
\vec{u}\cdot x+2t\partial_t \vec{u})^2-
(|\nabla \vec{u}|^2-2\chi_{\Omega}F(\vec{u}))\nabla \varphi \cdot xG]\nonumber\\
&\quad+r^{-2\beta-1}\int_{T^-_r}2[G\nabla \vec{u} \cdot \nabla
\varphi \nabla \vec{u}\cdot x+\nabla \vec{u}\cdot \nabla \varphi 2t
\partial_t \vec{u} G].
\end{align*}
In another word,
\begin{align} 0&=r^{-2\beta-1}[\int_{\mathbb{R}^n}2t\varphi
G(|\nabla
\vec{u}|^2-2\chi_{\Omega}F(\vec{u}))]^{(-r^2)}_{(-4r^2)}\label{3.10}\\
&\quad-2\beta r^{-2\beta-1}\int_{T^-_r}\varphi G(|\nabla
\vec{u}|^2-2\chi_{\Omega}F(\vec{u}))\nonumber\\
&\quad+r^{-2\beta-1}\int_{T^-_r}\frac{\varphi G}{t}(\nabla u\cdot
x+2t\partial_t \vec{u}-\beta \vec{u})^2\nonumber\\
&\quad +2r^{-2\beta-1}\int_{T^-_r}\chi_{\Omega}[\beta
\vec{u}\vec{f}(\vec{u})-2(\beta-1)F(\vec{u})]\varphi G\nonumber\\
&\quad+r^{-2\beta-1}\int_{T^-_r}\varphi
G(\frac{\beta}{t}\vec{u}\nabla
\vec{u}\cdot x+2\beta \vec{u}\partial_t \vec{u}-\frac{\beta^2}{t}\vec{u}^2)\nonumber\\
&\quad -r^{-2\beta-1}\int_{T^-_r}[2\beta \vec{u}G\nabla \varphi
\cdot
\nabla \vec{u}-(|\nabla \vec{u}|^2-2\chi_{\Omega}F(\vec{u}))\nabla \varphi \cdot xG]\nonumber\\
&\quad +r^{-2\beta-1}\int_{T^-_r}[2G\nabla \vec{u}\cdot \nabla
\varphi \nabla \vec{u}\cdot x+4tG\partial_t \vec{u}\nabla
\vec{u}\cdot \nabla \varphi].\nonumber
\end{align}
Meanwhile, we can see that
\begin{align}
&\quad r^{-2\beta-1}\int_{T^-_r}\varphi
G(\frac{\beta}{t}\vec{u}\nabla \vec{u}\cdot
x+2\beta \vec{u}\partial_t \vec{u}-\frac{\beta^2}{t}\vec{u}^2)\label{3.11}\\
&=\beta
\int_{T^-_1}\varphi(rx)\frac{G(t,x)}{t}r^{-\beta}\vec{u}(r^2t,rx)
[r^{-\beta}(\nabla \vec{u})(r^2t,rx)\cdot x\nonumber\\
&\quad +r^{-\beta+1}2t(\partial_t\vec{u})(r^2t,rx)-\beta
r^{-\beta-1}\vec{u}(r^2t,rx)]\nonumber\\
&=\partial_r(\frac{\beta}{2}\int_{T^-_1}(\frac{\vec{u}(r^2t,rx)}
{r^{\beta}})^2\frac{G(t,x)}{t}\varphi(rx))
-\frac{\beta}{2}r^{-2\beta-1}\int_{T^-_r}\frac{G\vec{u}^2}{t}\nabla \varphi\cdot x\nonumber\\
&=\partial_r(\frac{\beta}{2}\int_{T^-_r}(\frac{\vec{u}(x,t)}
{r^{\beta}})^2\frac{G(t,x)}{t}\varphi(x))
-\frac{\beta}{2}r^{-2\beta-1}\int_{T^-_r}\frac{G\vec{u}^2}{t}\nabla
\varphi\cdot x.\nonumber
\end{align}
Using (\ref{3.10}), (\ref{3.11}), we can obtain that
\begin{align*}
\frac{d}{dr}\Psi^-(r)&\geq 2r^{-2\beta-1}\int_{T^-_r}\chi_{\Omega}
[2(\beta-1)F(\vec{u})-\beta \vec{u}\vec{f}(\vec{u})]\varphi G\\
&\quad +r^{-2\beta-1}\int_{T^-_r}\frac{\varphi G}{(-t)}(\nabla
\vec{u}\cdot x+2t\partial_t \vec{u}-\beta \vec{u})^2+I,
\end{align*}
where
\begin{align*}
I&:=r^{-2\beta-1}\int_{T^-_r}[2\beta \vec{u}G\nabla \varphi \cdot
\nabla \vec{u}-(|\nabla \vec{u}|^2-2\chi_{\Omega}F(\vec{u}))\nabla \varphi \cdot xG]\nonumber\\
&\quad -r^{-2\beta-1}\int_{T^-_r}[2G\nabla \vec{u}\cdot \nabla
\varphi \nabla \vec{u}\cdot x+4tG\partial_t \vec{u}\nabla
\vec{u}\cdot \nabla
\varphi]\nonumber\\
&\quad +\frac{\beta}{2}r^{-2\beta-1}\int_{T^-_r}\frac{G\vec{u}^2}{t}
x\cdot \nabla \varphi.
\end{align*}

Since $u$ satisfies ($\star$) and $supp\varphi\subset B_{3/4}$, we
can see that the integrand in $I$ vanishes a.e. in $B_{1/2}\times
[-1,0]$ and $B^c_{3/4}\times [-1,0]$. Hence, we find that
\begin{align*}
I&\geq -r^{-2\beta-1}\int^{-r^2}_{-4r^2} \int_{B_{3/4}\setminus
B_{1/2}} (h_1(x,t)+\frac{h_2(x,t)}{t})G(x,-t)dxdt\\
&\geq-Cr^{-n-2\beta-1}e^{-1/(16r^2)}
\end{align*}
 with $\|h_1 \|_{L^1(Q^-_{3/4})}, \|h_2
\|_{L^1(Q^-_{3/4})}\leq C=C(n, \varphi, \beta)<\infty$ and
consequently
\begin{align*}
\frac{d}{dr}\Psi^-(r)&\geq
2r^{-2\beta-1}\int_{T^-_r}\chi_{\Omega}\varphi G
(2(\beta-1)F(\vec{u})-\beta \vec{u}\vec{f}(\vec{u}))\\
&\quad +r^{-2\beta-1}\int_{T^-_r}\frac{\varphi G}{(-t)}(\nabla
\vec{u}\cdot x+2t\partial_t \vec{u}-\beta \vec{u})^2\\
&\quad -Cr^{-n-2\beta-1}e^{-1/(16r^2)}.
\end{align*}
Therefore the function
$$\Psi^-(r)+CE(r)$$
is nondecreasing, where
$$
E(r)=\int^r_0 s^{-n-2\beta-1}e^{-1/(16s^2)}ds.$$

\end{proof}

We now study the blow-up of solutions. For a given point $(T,x_0)$
and a given sequence $\rho_k\rightarrow 0$, we define the {\it
scaled sequences} as follows:
$$
\vec{u}_k(t,x):=\rho^{-\beta}_k \vec{u}(T+\rho^2_kt,x_0+\rho_kx)
$$
and want to obtain more information on the solution's behavior. In
fact, we find that

\begin{theorem}
Suppose that for $t_1\leq T\leq t_2$ and $x_0\in \mathbb{R}^n$,
\begin{align*}
&\quad {\rm
sup}_{t\in(t_1,T-\delta)\cup(T+\delta,t_2)}\int_{\mathbb{R}^n}
\exp(-\frac{|x-x_0|^2}{4(T-t)})(|\nabla \vec{u}|^2-2F(\vec{u}))(t)\\
&+\int_{t\in(t_1,T-\delta)\cup(T+\delta,t_2)}\int_{\mathbb{R}^n}
\exp(-\frac{|x-x_0|^2}{4(T-t)})((\partial_t\vec{u})^2+\vec{u}^2)<\infty
\end{align*}
for any positive $\delta$, where $\vec{u}\in H^1(((t_1,T)\cup
(T,t_2))\times \mathbb{R}^n))$ is a variational solution
of(\ref{3.1}).

Suppose, furthermore, that in either case the growth estimates
\begin{align*}
&{\rm sup}_{r\in(0,\frac{\sqrt{T-t_1}}{4})}\max(r^{-2\beta}
\int_{T^-_r}\frac{1}{T-t}\vec{u}^2G_{(T,x_0)}\\
&\qquad +r^{-2\beta}\int_{T^-_r} |\nabla
\vec{u}|^2G_{(T,x_0)})+r^{-2\beta}|\int_{T^-_r}
2F(\vec{u}))G_{(T,x_0)}|<\infty
\end{align*}
and
\begin{align*}
&{\rm sup}_{r\in(0,\frac{\sqrt{t_2-T}}{4})}\max(r^{-2\beta}
\int_{T^+_r}\frac{1}{T-t}u^2G_{(T,x_0)}\\
&\qquad +r^{-2\beta}\int_{T^+_r} |\nabla
\vec{u}|^2G_{(T,x_0)}+r^{-2\beta}|\int_{T^+_r}
2F(\vec{u}))G_{(T,x_0)}|<\infty
\end{align*}
are satisfied. Then $\Psi^-(r)\searrow M^-(\vec{u},(T,x_0))$ as
$r\searrow 0$ provided that $T>t_1$ and $\Psi^+(r)\searrow
M^+(\vec{u},(T,x_0))$ as $r\searrow 0$ provided that $T<t_2$, and
for any $D\subset\subset
(((-\infty)\sqrt{T-t_1},0)\cup(0,((+\infty)\sqrt{t_2-T})))\times
\mathbb{R}^n$ and $k\geq k(D)$ the sequence
$$
\vec{u}_k(t,x):=\rho^{-\beta}_k \vec{u}(T+\rho^2_kt,x_0+\rho_kx)$$
is bounded in $H^1(D)\cap L^2(D)$ and any weak $H^1$-limit
$\vec{u}_0$ with respect to a subsequence is a function homogeneous
of degree $\beta$ on paths $\theta\rightarrow (\theta^2t,\theta x)$
for $\theta>0$ and $(t,x)\in
(((-\infty)\sqrt{T-t_1},0)\cup(0,((+\infty)\sqrt{t_2-T})))\times
\mathbb{R}^n$, i.e.,
$$
\vec{u}_0(\lambda^2t,\lambda x)=\lambda^{\beta}\vec{u}_0(t,x) \quad
{\rm for \ any} \ \lambda>0 $$ and for any $
(((-\infty)\sqrt{T-t_1},0)\cup(0,((+\infty)\sqrt{t_2-T})))\times
\mathbb{R}^n$.
\end{theorem}
\begin{proof}
We give the proof only for the case $t_2=T$ to avoid clumsy
notation. Calculating for $0<R<\infty$ that
\begin{align*}
\Psi^-(\rho_kR)&=R^{-2\beta}\int_{T^-_R(0)}|\nabla
\vec{u}_k|^2G_{(0,0)}-(\rho_k R)^{-2\beta}\int_{T^-_{\rho_k R}} 2F(\vec{u}))G_{(0,0)}\\
&\quad-\frac{\beta}{2}R^{-2\beta}
\int_{T^-_R(0)}\frac{1}{(-t)}\vec{u}^2_kG_{(0,0)},
\end{align*}
we know that the sequence $\vec{u}_k$ and $\nabla \vec{u}_k$ are
bounded in $L^2(D)$ for $k\geq k(D)$ by the assumed growth
estimate and the monotonicity formula Theorem 3.1.

By the results of Theorem 3.1, we know that $\Psi^-$ is
nondecreasing and bounded in $(0,r_0)$ for small positive $r_0$,
which means that $\Psi^-$ has a real right limit at $0$ and for
$0<R<S<\infty$,
\begin{align*}
0&\leftarrow \Psi^-(\rho_k S)-\Psi^-(\rho_k R)\\
&\quad =\int^{\rho_k S}_{\rho_k
R}2r^{-2\beta-1}\int_{T^-_r}(2(\beta-1)F(\vec{u})-\beta \vec{u}\vec{f}(\vec{u}))G_{(0,0)}\\
&\quad +\int^S_R r^{-2\beta-1}\int_{T^-_r}\frac{1}{(-t)}(\nabla
\vec{u}_k\cdot x+2t\partial_t \vec{u}_k-\beta \vec{u}_k)^2G_{(0,0)}.
\end{align*}
Then we can get that
\begin{align*}
0&\leftarrow \int^S_R r^{-2\beta-1}\int_{T^-_r}\frac{1}{(-t)}(\nabla
\vec{u}_k\cdot x+2t\partial_t \vec{u}_k-\beta \vec{u}_k)^2G_{(0,0)}
\end{align*}
as $k\rightarrow \infty$. Thus for $k\geq k(D)$ the sequence
$\vec{u}_k$ is bounded in $H^1(D)$. Since the lower
semi-continuity of the $L^2$-norm with respect to weak
convergence, we can take a subsequence $k\rightarrow \infty$ such
that $\vec{u}_k\rightharpoonup \vec{u_0}$ weakly convergence, and
obtain that
$$\nabla \vec{u}_0(t,x)\cdot x+2t\partial_t \vec{u}_0(t,x)-\beta
\vec{u}_0(t,x)=0$$
 a.e. in $(-\infty, 0)\times \mathbb{R}^n$. Now we
can easily see that $\vec{u_0}$ is homogeneous of degree $\beta$ on
paths $\theta\rightarrow (\theta^2t, \theta x)$ for $\theta>0$ and
$(t,x)\in (-\infty,0)\times \mathbb{R}^n$.
\end{proof}

\section{the monotonicity formulae for (\ref{0.3}) and (\ref{0.4}) }
In this section, we exhibit our monotonicity formulae for
 (\ref{0.3}) and (\ref{0.4}) in some special cases
 since they take simple forms. We will show that our monotonicity
 formulae do give new results even for free boundary problems
 related to those considered by Weiss \cite{W}.

 \subsection{The monotonicity formulae for (\ref{0.3})}

 First, we consider the elliptic equation case. Here we have to replace
 $\vec{u}$ by $u$, $\vec{f}$ by $f$ in the corresponding results of
 section 2. Then we can obtain the following theorems.
\begin{theorem}
Assume that $u$ is a solution of (\ref{0.3}) in the sense of
variations, $B_{\delta}(x_0)\subset\subset \Omega$ with $0<\delta$.
Then for any $\beta$ such that
 $\int_{B_r(x_0)}[2(\beta-1)F(u)-\beta
uf(u)]\geq 0$ and for all $ 0<\rho<\sigma<\delta$ the function $$
\Phi_{x_0}(r):=r^{-n-2\beta+2}\int_{B_r(x_0)} (|\nabla
u|^2-2F(u))-\beta r^{-n-2\beta+1}\int_{\partial B_r(x_0)}
u^2d\mathcal{H}^{n-1},
$$
defined in $(0,\delta)$, is nondecreasing in $r$ and satisfies the
monotonicity formula
\begin{align*}
\Phi_{x_0}(\sigma)-\Phi_{x_0}(\rho)&=
\int_{\rho}^{\sigma}2r^{-n-2\beta+1}\int_{B_r(x_0)}[2(\beta-1)F(u)-\beta
uf(u)]\\
&\ +\int_{\rho}^{\sigma}2r^{-n-2\beta+2}\int_{\partial B_r(x_0)}
(\nabla u\cdot \nu-\beta \frac{u}{r})^2 d\mathcal{H}^{n-1}dr\geq
0.\nonumber
\end{align*}
\end{theorem}

Now we give some examples to explain the Theorem.

\begin{corollary}
Assume that $u$ is a variational solution of
$$
\Delta u+u^p=0 \quad {\rm in} \ \Omega
$$
where $p\neq \pm 1$, and $B_{\delta}(x_0)\subset\subset
\Omega\subset \subset \mathbb{R}^n$. Then, for all $\beta$
satisfying $ \frac{\beta(1-p)-2}{(p+1)}\geq 0$ and for $
0<\rho<\sigma<\delta$, the function
$$
\Phi_{x_0}(r):=r^{-n-2\beta+2}\int_{B_r(x_0)} (|\nabla
u|^2-\frac{2u^{p+1}}{p+1})-\beta r^{-n-2\beta+1}\int_{\partial
B_r(x_0)} u^2d\mathcal{H}^{n-1},
$$
defined in $(0,\delta)$, is nondecreasing in $r$ and satisfies the
monotonicity formula
\begin{align*}
\Phi_{x_0}(\sigma)-\Phi_{x_0}(\rho)&=
\int_{\rho}^{\sigma}2r^{-n-2\beta+1}\int_{B_r(x_0)}(
\frac{2(\beta-1)}{(p+1)}-\beta)u^{p+1}\\
&\ +\int_{\rho}^{\sigma}2r^{-n-2(\beta-1)}\int_{\partial B_r(x_0)}
(\nabla u\cdot \nu-\beta \frac{u}{r})^2 d\mathcal{H}^{n-1}dr\geq
0.\nonumber
\end{align*}
\end{corollary}

{\bf Remark 4.1.} Since $f(u)=u^p$, $F(u)=\frac{u^{p+1}}{p+1}$ for
$p\neq \pm 1$, and $\frac{\beta(1-p)-2}{p+1}\geq 0$, we have
\begin{align*}
\quad \int_{B_r(x_0)}(2(\beta-1)F(u)-\beta
uf(u))=\int_{B_r(x_0)}\frac{\beta(1-p)-2}{p+1}u^{p+1}\geq 0.
\end{align*}

{\bf Remark 4.2.} The constant $\beta$ may be positive or negative
in the inequality of $\frac{\beta(1-p)-2}{p+1}\geq 0$. In fact, we
can choose $\beta\geq 2/(1-p)\geq 0$ if $-1<p<1$, and $\beta\leq
2/(1-p)\leq 0$ if $p>1$, and $\beta\leq 2/(1-p)$ being positive or
not if $p<-1$.

\begin{corollary}
Assume that $u$ is a variational solution of
$$
\Delta u+u^{-1}=0 \quad {\rm in} \ \Omega
$$
and $B_{\delta}(x_0)\subset\subset \Omega\subset \subset
\mathbb{R}^n$. Then, for $0<r<\delta$ and for any real $c$, we can
choose the constant $\beta$ such that
$$
\int_{B_r(x_0)}( 2(\beta-1)(\log u+c)-\beta)\geq 0.
$$
Therefore,
$$
\Phi_{x_0}(r)=r^{-n-2\beta+2}\int_{B_r(x_0)} (|\nabla u|^2-2\log
u-2c)-\beta r^{-n-2\beta+1}\int_{\partial B_r(x_0)}
u^2d\mathcal{H}^{n-1}
$$ is increasing in $r$ and satisfies
\begin{align*}
\Phi_{x_0}(\sigma)-\Phi_{x_0}(\rho)&=
\int_{\rho}^{\sigma}2r^{-n-2\beta+1}\int_{B_r(x_0)}(
2(\beta-1)(\log u+c)-\beta)\\
&\ +\int_{\rho}^{\sigma}2r^{-n-2(\beta-1)}\int_{\partial B_r(x_0)}
(\nabla u\cdot \nu-\beta \frac{u}{r})^2 d\mathcal{H}^{n-1}dr\geq
0.\nonumber
\end{align*}
\end{corollary}

We can also characterize the scaled sequences as follows:
\begin{theorem}
Suppose that $0<\rho_k \rightarrow 0$ as $k\rightarrow \infty$,
$u$ is in $B_{\delta}(x_0)$ a variational solution of (\ref{0.3}),
and that $u$ satisfies in $x_0$ the growth estimate
$$
{\rm sup}_{r\in (0,\delta)}(r^{-n-2\beta+1}\int_{\partial B_r(x_0)}
u^2d\mathcal{H}^{n-1}+r^{-n-2\beta+2}|\int_{B_r(x_0)}2F(u)|)<\infty.
$$
Then $\Phi_{x_0}(r)\searrow M(u,x_0)$ as $r\rightarrow 0$ and for
any open $D\subset \subset \mathbb{R}^n$ and $k\geq k(D)$, the
sequence $u_k(x):=\rho_k^{-\beta}u(x_0+\rho_kx)$ is bounded in
$H^{1}(D)$ and any weak $H^{1}$-limit with respect to a subsequence
$k\rightarrow \infty$ is homogeneous of degree $\beta$.
\end{theorem}

\subsection{ The monotonicity formula for (\ref{0.4})}

We now give the monotonicity formula for (\ref{0.4}). Here we have
to replace  $\vec{u}$ by $u$, $\vec{f}$ by $f$ in the
corresponding results of section 3. Note that the monotonicity
formula holds as soon as the solution continues to exist.

\begin{theorem}
(monotonicity formula) Assume that for $t_1\leq T\leq t_2$,
$x_0\in \mathbb{R}^n$, we have
\begin{align*}
&\quad {\rm sup}_{t\in (t_1,
T-\delta)\cup(T+\delta,t_2)}\int_{\mathbb{R}^n}\exp(-\frac{|x-x_0|^2}{4(T-t)})(|\nabla
u|^2-2F(u))(t,x)dx\\
&+\int_{(t_1,T-\delta)\cup(T+\delta,t_2)}\int_{\mathbb{R}^n}\exp(-\frac{|x-x_0|^2}{4(T-t)})
((\partial_t u)^2+u^2)(t,x)dxdt<\infty
\end{align*}
for any positive $\delta$, where
$$u\in H^1(((t_1,T)\cup
(T,t_2))\times \mathbb{R}^n)$$ is a variational solution of
(\ref{0.4}). Then for any $\beta$ satisfying
$$\int_{T^-_r} G_{(T,x_0)}[2(\beta-1)F(u)-\beta uf(u)]\geq 0$$
and $$\int_{T^+_r} G_{(T,x_0)}[2(\beta-1)F(u)-\beta uf(u)]\geq 0,$$
the functions
\begin{align*}
\Psi^-(r)&:=r^{-2\beta}\int_{T^-_r(T)}(|\nabla
u|^2-2F(u))G_{(T,x_0)}\\
&\quad-\frac{\beta}{2}r^{-2\beta}\int_{T^-_r(T)}\frac{1}{T-t}u^2G_{(T,x_0)}
\end{align*}
and
\begin{align*}
\Psi^+(r)&:=r^{-2\beta}\int_{T^+_r(T)}(|\nabla
u|^2-2F(u))G_{(T,x_0)}\\
&\quad-\frac{\beta}{2}r^{-2\beta}\int_{T^+_r(T)}\frac{1}{T-t}u^2G_{(T,x_0)}
\end{align*}
are well defined in the interval $(0,\frac{\sqrt{T-t_1}}{2})$ and
$(0,\frac{\sqrt{t_2-T}}{2})$, respectively, and satisfy for any
$0<\rho<\sigma<\frac{\sqrt{T-t_1}}{2}$ and
$0<\rho<\sigma<\frac{\sqrt{t_2-T}}{2}$, respectively, the
monotonicity formulae
\begin{align*}
&\Psi^-(\sigma)-\Psi^-(\rho)=\int^{\sigma}_{\rho}2r^{-2\beta-1}
\int_{T^-_r(T)}[2(\beta-1)F(u)-\beta uf(u)]G_{(T,x_0)}\\
&+\int^{\sigma}_{\rho}r^{-2\beta-1}\int_{T^-_r(T)}\frac{1}{T-t}(\nabla
u\cdot(x-x_0)-2(T-t)\partial_t u-\beta u)^2G_{(T,x_0)}\\
&\geq 0
\end{align*}
and
\begin{align*}
&\Psi^+(\sigma)-\Psi^+(\rho)=\int^{\sigma}_{\rho}2r^{-2\beta-1}
\int_{T^+_r(T)}[2(\beta-1)F(u)-\beta uf(u)]G_{(T,x_0)}\\
&+\int^{\sigma}_{\rho}r^{-2\beta-1}\int_{T^+_r(T)}\frac{1}{T-t}(\nabla
u\cdot(x-x_0)-2(T-t)\partial_t u-\beta u)^2G_{(T,x_0)}\\
&\geq 0.
\end{align*}
\end{theorem}

We now turn our attention to other kinds of problems. Consider the
following free boundary problem:
$$
(*) \left\{ \begin{array}{ll} u_t-\Delta u=\chi_{\Omega} f(u) &
{\rm
in} \ \mathbb{R}^n\times \mathbb{R},\\
u=|\nabla u|=0 & {\rm in} \; \Lambda:=(\mathbb{R}^n\times
\mathbb{R}) \setminus \Omega.
\end{array} \right.
$$

We have that

\begin{theorem}
Assume that $u$ is a variational solution of $(*)$, $(x_0,T)\in
\Lambda$, $t_1\leq T\leq t_2$, and
\begin{align*}
&\quad {\rm sup}_{t\in (t_1,
T-\delta)\cup(T+\delta,t_2)}\int_{\mathbb{R}^n}\exp(-\frac{|x-x_0|^2}{4(T-t)})(|\nabla
u|^2-2\chi_{\Omega}F(u))(t,x)dx\\
&+\int_{(t_1,T-\delta)\cup(T+\delta,t_2)}\int_{\mathbb{R}^n}\exp(-\frac{|x-x_0|^2}{4(T-t)})
((\partial_t u)^2+u^2)(t,x)dxdt<\infty
\end{align*}
for any positive $\delta$. Let $\varphi(x)\geq 0$ be a $C^{\infty}$
cut-off function in $\mathbb{R}^n$ with ${\rm supp}\varphi\subset
B_{3/4}(x_0)$ and $\varphi|_{B_{1/2}}=1$. Then for any $\beta$
satisfying
$$\int_{T^-_r} \chi_{\Omega}[2(\beta-1)F(u)-\beta uf(u)]\varphi G_{(T,x_0)}\geq 0$$
and
$$\int_{T^+_r} \chi_{\Omega}[2(\beta-1)F(u)-\beta uf(u)]\varphi G_{(T,x_0)}\geq 0,$$
there exists a constant $C=C(n,\varphi,\beta)>0$ such that the
functions
\begin{align*}
\Psi^-(r)&:=r^{-2\beta}\int_{T^-_r(T)}(|\nabla
u|^2-2\chi_{\Omega}F(u))\varphi G_{(T,x_0)}\\
&\quad-\frac{\beta}{2}r^{-2\beta}\int_{T^-_r(T)}\frac{1}{T-t}u^2G_{(T,x_0)}\\
&\quad +C\int^r_0 s^{-n-2\beta-1}e^{-1/(16s^2)}ds
\end{align*}
and
\begin{align*}
\Psi^+(r)&:=r^{-2\beta}\int_{T^+_r(T)}(|\nabla
u|^2-2\chi_{\Omega}F(u)) \varphi G_{(T,x_0)}\\
&\quad-\frac{\beta}{2}r^{-2\beta}\int_{T^+_r(T)}\frac{1}{T-t}u^2\varphi
G_{(T,x_0)}\\
&\quad +C\int^r_0 s^{-n-2\beta-1}e^{-1/(16s^2)}ds
\end{align*}
are well defined in the interval $(0,\frac{\sqrt{T-t_1}}{2})$ and
$(0,\frac{\sqrt{t_2-T}}{2})$, respectively, and satisfy for any
$0<\rho<\sigma<\frac{\sqrt{T-t_1}}{2}$ and
$0<\rho<\sigma<\frac{\sqrt{t_2-T}}{2}$, respectively, the
monotonicity formulae
\begin{align*}
&\Psi^-(\sigma)-\Psi^-(\rho)=\int^{\sigma}_{\rho}2r^{-2\beta-1}
\int_{T^-_r(T)}\chi_{\Omega}[2(\beta-1)F(u)-\beta uf(u)]\varphi G_{(T,x_0)}\\
&+\int^{\sigma}_{\rho}r^{-2\beta-1}\int_{T^-_r(T)}\frac{1}{T-t}(\nabla
u\cdot(x-x_0)-2(T-t)\partial_t u-\beta u)^2\varphi G_{(T,x_0)}\\
&+C\int^{\sigma}_{\rho} s^{-n-2\beta-1}e^{-1/(16s^2)}ds\geq 0
\end{align*}
and
\begin{align*}
&\Psi^+(\sigma)-\Psi^+(\rho)=\int^{\sigma}_{\rho}2r^{-2\beta-1}
\int_{T^+_r(T)}\chi_{\Omega}[2(\beta-1)F(u)-\beta uf(u)]\varphi G_{(T,x_0)}\\
&+\int^{\sigma}_{\rho}r^{-2\beta-1}\int_{T^+_r(T)}\frac{1}{T-t}(\nabla
u\cdot(x-x_0)-2(T-t)\partial_t u-\beta u)^2\varphi G_{(T,x_0)}\\
&+C\int^{\sigma}_{\rho} s^{-n-2\beta-1}e^{-1/(16s^2)}ds\geq 0.
\end{align*}
\end{theorem}

We can give the characterizing of the {\it scaled sequence}. We have
that

\begin{theorem} Let
$$u\in H^1((t_1,T)\cup
(T,t_2))\times \mathbb{R}^n)$$ be in  a variational solution of
(\ref{0.4}). Suppose that for $t_1\leq T\leq t_2$ and $x_0\in
\mathbb{R}^n$, we have
\begin{align*}
&\quad {\rm
sup}_{t\in(t_1,T-\delta)\cup(T+\delta,t_2)}\int_{\mathbb{R}^n}
\exp(-\frac{|x-x_0|^2}{4(T-t)})(|\nabla u|^2-2F(u))(t)\\
&+\int_{t\in(t_1,T-\delta)\cup(T+\delta,t_2)}\int_{\mathbb{R}^n}
\exp(-\frac{|x-x_0|^2}{4(T-t)})((\partial_tu)^2+u^2)<\infty
\end{align*}
for any small positive $\delta$.

Suppose, furthermore, that in either case the growth estimates
\begin{align*}
&{\rm sup}_{r\in(0,\frac{\sqrt{T-t_1}}{4})}\max(r^{-2\beta}
\int_{T^-_r}\frac{1}{T-t}u^2G_{(T,x_0)}\\
&\qquad+r^{-2\beta}\int_{T^-_r} |\nabla
u|^2G_{(T,x_0)}+r^{-2\beta}|\int_{T^-_r} 2F(u)G_{(T,x_0)}|)<\infty
\end{align*}
and
\begin{align*}
&{\rm sup}_{r\in(0,\frac{\sqrt{t_2-T}}{4})}\max(r^{-2\beta}
\int_{T^+_r}\frac{1}{T-t}u^2G_{(T,x_0)}\\
&\qquad+r^{-2\beta}\int_{T^+_r} |\nabla
u|^2G_{(T,x_0)}+r^{-2\beta}|\int_{T^+_r} 2F(u))G_{(T,x_0)}|)<\infty
\end{align*}
are satisfied. Then $\Psi^-(r)\searrow M^-(u,(T,x_0))$ as
$r\searrow 0$ provided that $T>t_1$ and $\Psi^+(r)\searrow
M^+(u,(T,x_0))$ as $r\searrow 0$ provided that $T<t_2$, and for
any $D\subset\subset
(((-\infty\sqrt{T-t_1},0)\cup(0,((+\infty)\sqrt{t_2-T})))\times
\mathbb{R}^n$ and $k\geq k(D)$ the sequence
$$
u_k(t,x):=\rho^{-\beta}_k u(T+\rho^2_kt,x_0+\rho_kx)$$ is bounded in
$H^1(D)\cap L^2(D)$ and any weak $H^1$-limit $u_0$ with respect to a
subsequence is a function homogeneous of degree $\beta$ on paths
$\theta\rightarrow (\theta^2t,\theta x)$ for $\theta>0$ and
$$(t,x)\in
(((-\infty)\sqrt{T-t_1},0)\cup(0,((+\infty)\sqrt{t_2-T})))\times
\mathbb{R}^n,$$ i.e.,
$$
u_0(\lambda^2t,\lambda x)=\lambda^{\beta}u_0(t,x) \quad {\rm for \
any} \ \lambda>0 $$ and for any $(x,t)\in
(((-\infty)\sqrt{T-t_1},0)\cup(0,((+\infty)\sqrt{t_2-T})))\times
\mathbb{R}^n$.
\end{theorem}

We now give some examples to explain these Theorems.

\begin{corollary}
Assume that $t_1<T<t_2$, $x_0\in \mathbb{R}^n$, and
$$u\in H^1(((t_1,T)\cup (T,t_2))\times \mathbb{R}^n)$$ is a variational
solution of
$$u_t-\Delta u=u^p \quad {\rm in} \ (t_1, t_2)\times
\mathbb{R}^n,$$ where $p \neq \pm 1$. Suppose furthermore that
\begin{align*}
&\quad {\rm sup}_{t\in (t_1,
T-\delta)\cup(T+\delta,t_2)}\int_{\mathbb{R}^n}\exp(-\frac{|x-x_0|^2}{4(T-t)})(|\nabla
u|^2-\frac{2u^{p+1}}{p+1})dx\\
&+\int_{(t_1,T-\delta)\cup(T+\delta,t_2)}\int_{\mathbb{R}^n}\exp(-\frac{|x-x_0|^2}{4(T-t)})
((\partial_t u)^2+u^2)(t,x)dxdt<\infty
\end{align*}
for any positive $\delta$, and $\frac{\beta(1-p)-2}{p+1}\geq 0$.
Then the functions
\begin{align*}
\Psi^-(r)&:=r^{-2\beta}\int_{T^-_r(T)}(|\nabla
u|^2-\frac{2u^{p+1}}{p+1})G_{(T,x_0)}\\
&\quad-\frac{\beta}{2}r^{-2\beta}\int_{T^-_r(T)}\frac{1}{T-t}u^2G_{(T,x_0)}
\end{align*}
and
\begin{align*}
\Psi^+(r)&:=r^{-2\beta}\int_{T^+_r(T)}(|\nabla
u|^2-\frac{2u^{p+1}}{p+1})G_{(T,x_0)}\\
&\quad-\frac{\beta}{2}r^{-2\beta}\int_{T^+_r(T)}\frac{1}{T-t}u^2G_{(T,x_0)}
\end{align*}
are well defined in the interval $(0,\frac{\sqrt{T-t_1}}{2})$ and
$(0,\frac{\sqrt{t_2-T}}{2})$, respectively, and satisfy for any
$0<\rho<\sigma<\frac{\sqrt{T-t_1}}{2}$ and
$0<\rho<\sigma<\frac{\sqrt{t_2-T}}{2}$, respectively, the
monotonicity formulae
\begin{align*}
&\Psi^-(\sigma)-\Psi^-(\rho)=\int^{\sigma}_{\rho}2r^{-2\beta-1}
\int_{T^-_r(T)}\frac{\beta(1-p)-2}{p+1}u^{p+1}\\
&+\int^{\sigma}_{\rho}r^{-2\beta-1}\int_{T^-_r(T)}\frac{1}{T-t}(\nabla
u\cdot(x-x_0)-2(T-t)\partial_t u-\beta u)^2G_{(T,x_0)}\\
&\geq 0
\end{align*}
and
\begin{align*}
&\Psi^+(\sigma)-\Psi^+(\rho)=\int^{\sigma}_{\rho}2r^{-2\beta-1}
\int_{T^+_r(T)}\frac{\beta(1-p)-2}{p+1}u^{p+1}\\
&+\int^{\sigma}_{\rho}r^{-2\beta-1}\int_{T^+_r(T)}\frac{1}{T-t}(\nabla
u\cdot(x-x_0)-2(T-t)\partial_t u-\beta u)^2G_{(T,x_0)}\\
&\geq 0.
\end{align*}
\end{corollary}

\begin{corollary}
Assume that $t_1<T<t_2$, $x_0\in \mathbb{R}^n$, and
$$u\in H^1(
((t_1,T)\cup (T,t_2))\times \mathbb{R}^n)$$ is a variational
solution of
$$u_t-\Delta u=u^{-1} \quad {\rm in} \ (t_1, t_2)\times
\mathbb{R}^n.$$ Suppose furthermore that
\begin{align*}
&\quad {\rm sup}_{t\in (t_1,
T-\delta)\cup(T+\delta,t_2)}\int_{\mathbb{R}^n}\exp(-\frac{|x-x_0|^2}{4(T-t)})(|\nabla
u|^2-2\log u)dx\\
&+\int_{(t_1,T-\delta)\cup(T+\delta,t_2)}\int_{\mathbb{R}^n}\exp(-\frac{|x-x_0|^2}{4(T-t)})
((\partial_t u)^2+u^2)(t,x)dxdt<\infty
\end{align*}
for any positive $\delta$. For any constant $c$, we can choose
$\beta$ such that
$$
\int_{T^-_r}( 2(\beta-1)(\log u+c)-\beta)\geq 0
$$
and
$$
\int_{T^+_r}( 2(\beta-1)(\log u+c)-\beta)\geq 0.
$$
Then the functions
\begin{align*}
\Psi^-(r)&:=r^{-2\beta}\int_{T^-_r(T)}(|\nabla
u|^2-2\log u-2c)G_{(T,x_0)}\\
&\quad-\frac{\beta}{2}r^{-2\beta}\int_{T^-_r(T)}\frac{1}{T-t}u^2G_{(T,x_0)}
\end{align*}
and
\begin{align*}
\Psi^+(r)&:=r^{-2\beta}\int_{T^+_r(T)}(|\nabla
u|^2-2\log u-2c)G_{(T,x_0)}\\
&\quad-\frac{\beta}{2}r^{-2\beta}\int_{T^+_r(T)}\frac{1}{T-t}u^2G_{(T,x_0)}
\end{align*}
are well defined in the interval $(0,\frac{\sqrt{T-t_1}}{2})$ and
$(0,\frac{\sqrt{t_2-T}}{2})$, respectively, and satisfy for any
$0<\rho<\sigma<\frac{\sqrt{T-t_1}}{2}$ and
$0<\rho<\sigma<\frac{\sqrt{t_2-T}}{2}$, respectively, the
monotonicity formulae
\begin{align*}
&\Psi^-(\sigma)-\Psi^-(\rho)=\int^{\sigma}_{\rho}2r^{-2\beta-1}
\int_{T^-_r(T)}(2(\beta-1)(\log u+c)-\beta)\\
&+\int^{\sigma}_{\rho}r^{-2\beta-1}\int_{T^-_r(T)}\frac{1}{T-t}(\nabla
u\cdot(x-x_0)-2(T-t)\partial_t u-\beta u)^2G_{(T,x_0)}\\
&\geq 0
\end{align*}
and
\begin{align*}
&\Psi^+(\sigma)-\Psi^+(\rho)=\int^{\sigma}_{\rho}2r^{-2\beta-1}
\int_{T^+_r(T)}(2(\beta-1)(\log u+c)-\beta)\\
&+\int^{\sigma}_{\rho}r^{-2\beta-1}\int_{T^+_r(T)}\frac{1}{T-t}(\nabla
u\cdot(x-x_0)-2(T-t)\partial_t u-\beta u)^2G_{(T,x_0)}\\
&\geq 0.
\end{align*}

\end{corollary}

So much for examples. It is clear from our examples that one can
use the monotonicity formulae in different forms for different
purposes. We hope to see more applications of them in the future.


\begin{thebibliography}{1}
\bibitem{Al}
W. K. Allard, \emph{On the first variation of a varifold}, Ann. of
Math. 95(1972), 417--491.

\bibitem{ACF}
H. W. Alt, L. A. Caffarelli, A. Friedman, \emph{Variational problems
with two phases and their free boundaries}, Trans. Amer. Math. Soc.
282(1984), 431--462.

\bibitem{BBM}
Bourgain, J.; Brezis, H.; Mironescu, P. \emph{$H^{1/2}$ maps with
values into the circle: minimal connections, lifting, and the
Ginzburg-Landau equation}. Publications math谷matiques de l＊IHES 89
(2004), 1每115.

\bibitem{Ca}
L. A. Caffarelli, \emph{A monotonicity formula for heat functions in
disjoint domains}, Boundary value problems for partial differential
equations and applications, Masson, Paris, 1993, pp. 53-60.

\bibitem{CJK}
L. Caffarelli, D. Jerison and C. E. Keng, \emph{Some new
monotonicity theorems with appliacations to free boundary problems},
Annals of Math. 155(2002), 369--402.

\bibitem{CK}
L. A. Caffarelli and C. E. Keng, \emph{Gradient estimates for
variable coefficient parabolic equations and singular perturbation
problems}, Amer. J. of Math. 120(1998), 391--439.

\bibitem{CKS}
L. Caffarelli, L. Karp and H. Shahgholian, \emph{Regularity of a
free boundary with application to the Pompeiu problem}, Annals of
Math. 151(2000), 269--292.

\bibitem{CPS}
L. Caffarelli, A. Petrosyan and H. Shahgholian, \emph{Regularity of
a free boundary in parabolic potential theory}. J. Amer. Math. Soc.
17 (2004), 827--869.

\bibitem{E1}
K. Ecker, \emph{A local monotonicity formula for mean curvature
flow}, Annals of Math. 154(2001), 503--523.

\bibitem{E2}
K. Ecker, \emph{Local monotonicity formulas for some nonlinear
diffusion equations}. Calc. Var. Partial Differential Equations 23
(2005), no. 1, 67--81.

\bibitem{Fl}
W. H. Fleming, \emph{On the oriented Plateau problem}, Rend. Circ.
Mat. Palermo(2)11(1962), 69--90.

\bibitem{GK}
Y. Giga and R. V. Kohn, \emph{Asymptotically self-similar blow-up of
semilinear heat equations}, Comm. Pure Appl. Math. 38(1985),
297--319.



\bibitem{Ha}
R. Hamilton, \emph{Monotonicity formulas for parabolic flows on
manifolds}, Comm. Analysis Geom. 1(1993), 127--137.

\bibitem{Hu}
G. Huisken, \emph{Asymptotic behaviour for singularities of the mean
curvature flow}, J. Differ. Geom. 31(1990), 285--299.

\bibitem{LSU}
O. A. Ladyzenskaja, V. A. Solonnikov and N. N. Ural'ceva,
\emph{Linear and Quasi-linear Equations of Parabolic Type}, Transl.
Math. Monogr. 23, American Mathematical Society, Providence, RI,
1988.

\bibitem{LR}
F. H. Lin and T. Riviere, \emph{Complex Ginzburg每Landau equations
in high dimensions and codimension two area minimizing currents}, J.
Eur. Math. Soc., 1 (1999), 237每311; Erratum 2 (2002), 87每91.

\bibitem{Pa}
F. Pacard, \emph{Partial regularity for weak solutions of a
nonlinear elliptic equation}, Manuscripta Math. 79(1993), 161-172.

\bibitem{Pr}
P. Price, \emph{A monotonicity formula for Yang-Mills fields},
Manuscripta Math. 43(1983), 131--166.

\bibitem{RI}
T. Riviere, \emph{Line vortices in the U(1)-Higgs model}, Control
Optim. Calc. Var., 1 (1996), 77每167.

\bibitem{Si}
L. M. Simon, \emph{Lectures on Geometric Measure Theory}, Proc. of
the CMA, Vol.3(1983).

\bibitem{St}
M. Struwe, \emph{On the evolution of harmonic maps in higher
dimensions}, J. Differ. Geom. 28(1988), 485--502.

\bibitem{SU}
R. M. Schoen, \emph{Analytic aspects of the harmonic map problem,}
Seminar on Nonlinear PDE, Springer, Berlin(1984).

\bibitem{W}
G. S. Weiss, \emph{Partial regularity for a minimum problem with
free boundary}. (English. English summary) J. Geom. Anal. 9 (1999),
no. 2, 317--326.

\bibitem{W1}
G. S. Weiss, \emph{Partial regularity for weak solution of an
elliptic free boundary problem}, Commun. Partial Differ. Equations
23(1998), 439--453.

\bibitem{W2} G. S. Weiss, \emph{A homogeneity improvement approach to the
obstacle problem}, Invent. Math. 138(1999), 23--50.

\bibitem{W3}
G. S. Weiss, \emph{Self-similar blow-up and Hausdorff dimension
estimates for a class of parabolic free boundary problems}, SIAM J.
Math. Anal. 30(1999), 623--642.

\end{thebibliography}
\end{document}